\documentclass[12pt]{article}
\usepackage{color}
\definecolor{darkblue}{rgb}{0.00,0.25,0.50}
\usepackage[colorlinks,filecolor=blue,citecolor=darkblue]{hyperref}

\usepackage{amsfonts,amssymb,amsmath,amsthm}
\usepackage{url}
\usepackage{enumerate}
\usepackage[ukrainian, russian, english]{babel}
\usepackage[cp1251]{inputenc}
\begin{document}\selectlanguage{english}
\thispagestyle{empty}

\title{}

\begin{center}
\textbf{\Large Uniform approximations by Fourier sums on  classes of generalized Poisson integrals }
\end{center}
\vskip0.5cm
\begin{center}
A.~S.~Serdyuk${}^1$, T.~A.~Stepanyuk${}^2$\\ \emph{\small
${}^1$Institute of Mathematics NAS of
Ukraine, Kyiv\\
${}^2$Lesya Ukrainka Eastern European National University, Lutsk\\}
\end{center}
\vskip0.5cm

%\address{Institute of Mathematics of the National Academy of
%Sciences of Ukraine\\ 3\\ Tereshenkivska st.\\ 01601\\ Kiev, Ukraine}

\begin{abstract}

We find asymptotic equalities for exact upper bounds of approximations by Fourier  sums in  uniform
metric on classes of
  $2\pi$--periodic functions,  representable in the form of convolutions of functions $\varphi$, which belong to unit balls of spaces $L_{p}$, with generalized Poisson kernels. For obtained asymptotic equalities we introduce the estimates of remainder, which are expressed in the explicit form via the parameters of the problem.
  
  \vskip 0.5cm
  

\noindent {\bf Math Subject Classifications} 42A05, 42A10, 42A16.
 
\end{abstract}

\vskip 0.5cm

{\bf 1. Introduction.}
Let $L_{p}$,
$1\leq p<\infty$, be the space of $2\pi$--periodic functions $f$ summable to the power $p$
on  $[0,2\pi)$, in which
the norm is given by the formula
${\|f\|_{p}=\Big(\int\limits_{0}^{2\pi}|f(t)|^{p}dt\Big)^{\frac{1}{p}}}$; $L_{\infty}$ be the space of measurable and essentially bounded   $2\pi$--periodic functions  $f$ with the norm
$\|f\|_{\infty}=\mathop{\rm{ess}\sup}\limits_{t}|f(t)|$; $C$ be the space of continuous $2\pi$--periodic functions  $f$, in which the norm is specified by the equality
 ${\|f\|_{C}=\max\limits_{t}|f(t)|}$.

Denote by $C^{\alpha,r}_{\beta,p}, \ \alpha>0, \ r>0, \ 1\leq p\leq\infty,$ the set of all  $2\pi$--periodic functions, representable for  all
$x\in\mathbb{R}$ as convolutions of the form (see, e.g., \cite[p.~133]{Stepanets1})
\begin{equation}\label{conv}
f(x)=\frac{a_{0}}{2}+\frac{1}{\pi}\int\limits_{-\pi}^{\pi}P_{\alpha,r,\beta}(x-t)\varphi(t)dt,
\ a_{0}\in\mathbb{R}, \ \varphi\in B_{p}^{0}, \
\end{equation}
$$
B_{p}^{0}=\left\{\varphi: \ ||\varphi||_{p}\leq 1, \  \varphi\perp1\right\},
 \ 1\leq p\leq \infty,
$$
with fixed generated kernels
$$
P_{\alpha,r,\beta}(t)=\sum\limits_{k=1}^{\infty}e^{-\alpha k^{r}}\cos
\big(kt-\frac{\beta\pi}{2}\big), \ \  \beta\in
    \mathbb{R}.
$$
The kernels  $P_{\alpha,r,\beta}(t)$ are called generalized Poisson kernels. For $r=1$ and $\beta=0$ the kernels $P_{\alpha,r,\beta}(t)$ are usual Poisson kernels of harmonic functions.

% (див., наприклад, [\ref{Step monog 1987}, с. 97])
For any $r>0$ the classes  $C^{\alpha,r}_{\beta,p}$ belong to set of infinitely differentiable
 $2\pi$--periodic functions $D^{\infty}$, i.e., $C^{\alpha,r}_{\beta,p}\subset D^{\infty}$ (see, e.g., \cite[p.~128]{Stepanets1}, \cite{Stepanets_Serdyuk_Shydlich}).
For  $r\geq1$ the classes  $C^{\alpha,r}_{\beta,p}$
consist of functions  $f$,  admitting a regular extension into the strip $|\mathrm{Im} \ z|\leq c, \ c>0$ in the complex
plane (see, e.g., \cite[p.~141]{Stepanets1}), i.e., are the classes of analytic functions.
For $r>1$
the classes  $C^{\alpha,r}_{\beta,p}$  consist of functions  regular on the whole complex plane,
i.e., of entire functions (see, e.g., \cite[p.~131]{Stepanets1}). Besides,  it follows from the theorem 1 in \cite{Stepanets_Serdyuk_Shydlich2009} that for any $r>0$ the embedding holds  $C^{\alpha,r}_{\beta,p}\subset \mathcal{J}_{1/r}$, where $\mathcal{J}_{a}, a>0,$ are known Gevrey classes
$$
\mathcal{J}_{a}=\bigg\{f\in D^{\infty}: \ \sup\limits_{k\in \mathbb{N}}\Big(\frac{\|f^{(k)}\|_{C}}{(k!)^{a}}\Big)^{1/k}<\infty \bigg\}.
$$

Approximation properties of classes of generalized Poisson integrals  $C^{\alpha,r}_{\beta,p}$ in metrics of spaces $L_{s}$, $1\leq s\leq \infty$,   were considered in    \cite{Kuchpel}--\cite{S_S2} from the  viewpoint of order or asymptotic estimates for approximations by Fourier sums, best approximations and widths.

In the present paper we obtain asymptotic equalities as  ${n\rightarrow\infty}$ for the quantities
  \begin{equation}\label{sum}
 {\cal E}_{n}(C^{\alpha,r}_{\beta,p})_{C}=\sup\limits_{f\in
C^{\alpha,r}_{\beta,p}}\|f(\cdot)-S_{n-1}(f;\cdot)\|_{C},  \ r>0, \ \alpha>0, \ 1\leq p \leq \infty,
  \end{equation}
  where $S_{n-1}(f;\cdot)$ are the partial Fourier sums of order  $n-1$ for a function $f$.

Approximation by Fourier sums on other classes of differentiable functions in uniform metric were investigated in works \cite{Stepanets1}, \cite{Kol}--\cite{Serdyuk_Stepaniuk2015}.

  Nikol’skii \cite[p.~221]{Nikolsky 1946} considered the case ${r=1}$, $p=\infty$ and established that following asymptotic equality is true
 \begin{equation}\label{triang}
  {\cal E}_{n}(C^{\alpha,1}_{\beta,\infty})_{C}=e^{-\alpha n}\Big(\frac{8}{\pi^{2}}{\bf K}(e^{-\alpha})+O(1)n^{-1}\Big),
 \end{equation}
      where
  $$
\mathbf{K}(q):=\int\limits_{0}^{\frac{\pi}{2}}\frac{dt}{\sqrt{1-q^{2}\sin^{2}t}}, \ q\in(0,1),
$$
is a complete elliptic integral of the first kind, and  $O(1)$ is a quantity uniformly bounded in parameters  $n$ and $\beta$.

  Later, the equality (\ref{triang}) was clarified by  Stechkin \cite[p.~139]{Stechkin 1980}, who established the asymptotic formula
\begin{equation}\label{stechkin}
  {\cal E}_{n}(C^{\alpha,1}_{\beta,\infty})_{C}=
 e^{-\alpha n}\Big(\frac{8}{\pi^{2}}\mathbf{K}(e^{-\alpha})+O(1)\frac{e^{-\alpha}}{(1-e^{-\alpha})n}\Big), \  \ \alpha>0, \ \beta\in\mathbb{R},
  \end{equation}
where  $O(1)$   is a quantity uniformly bounded  in all analyzed parameters.

In work \cite{Serdyuk2005} for $r=1$ and arbitrary values of  $1\leq p\leq\infty$ for quantities ${\cal E}_{n}(C^{\alpha,r}_{\beta,p})_{C}$, $\alpha>0$, $\beta\in\mathbb{R}$,
the following equality was established
\begin{equation}\label{hhj}
  {\cal E}_{n}(C^{\alpha,1}_{\beta,p})_{C}=e^{-\alpha n}\bigg(\frac{2}{\pi^{1+\frac{1}{p'}}}\|\cos t\|_{p'}K(p',e^{-\alpha})+O(1)\frac{e^{-\alpha}}{n(1-e^{-\alpha})^{s(p)}}\bigg), \ \ \end{equation}
where $p'=\frac{p}{p-1}$,
$$
s(p):={\left\{\begin{array}{cc}
1, \ & p=\infty,  \\
2, \ & p\in[1,2)\cup(2,\infty), \\
-\infty, &
p=2, \
  \end{array} \right.}
$$
$$
K(p',q):=\frac{1}{2^{1+\frac{1}{p'}}}\Big\|(1-2q\cos t+q^{2})^{-\frac{1}{2}} \Big\|_{p'}, \ q\in(0,1),
$$
and $O(1)$ is a quantity uniformly bounded in $n$, $p$, $\alpha$ and $\beta$.
For $p=\infty$, by virtue of  the known equality  $K(1,q)=\mathbf{K}(q)$, the estimate (\ref{hhj}) coincides with the estimate   (\ref{stechkin}).

Note that for $p=2$ and $r=1$ formula  (\ref{hhj}) becomes the equality
$$
 {\cal E}_{n}(C^{\alpha,1}_{\beta,2})_{C}=\frac{1}{\sqrt{\pi(1-e^{-2\alpha})}}e^{-\alpha n}, \ \alpha>0, \ \beta\in\mathbb{R}, \ n\in\mathbb{N},
$$
(see \cite{Serdyuk2005}).
Moreover,  it follows from  \cite{Serdyuk2011} that for $p=2$ and  $r>0$ for the quantities
${\cal E}_{n}(C^{\alpha, r}_{\beta,p})_{C}$ the equalities take place
\begin{equation}\label{serd2011}
 {\cal E}_{n}(C^{\alpha,r}_{\beta,2})_{C}=\frac{1}{\sqrt{\pi}}\Big(\sum\limits_{k=n}^{\infty}e^{-2\alpha k^{r}}\Big)^{\frac{1}{2}}, \ \alpha>0,  \ \beta\in\mathbb{R}, \ n\in\mathbb{N}.
\end{equation}

 In the case of ${r>1}$ and   $p=\infty$ the asymptotic equalities for the quantities  ${\cal E}_{n}(C^{\alpha,r}_{\beta,p})_{C}$, $\alpha>0$, $\beta\in\mathbb{R}$,   were obtained by Stepanets
 \cite[Chapter 3, Section~9]{Step monog 1987}, who showed that for any  $n\in\mathbb{N}$
\begin{equation}\label{step_prepr}
  {\cal E}_{n}(C^{\alpha,r}_{\beta,\infty})_{C}=\Big(\frac{4}{\pi}+\gamma_{n}\Big)e^{-\alpha n^{r}},
\end{equation}
where
$$
|\gamma_{n}|<2\Big(1+\frac{1}{\alpha r n^{r-1}}\Big)e^{-\alpha r n^{r-1}}.
$$

 Later Telyakovskii  \cite{Teljakovsky 1989} established the asymptotic equality
\begin{equation}\label{tel}
  {\cal E}_{n}(C^{\alpha,r}_{\beta,\infty})_{C}=\frac{4}{\pi}e^{-\alpha n^{r}}+
  O(1)\Big(e^{-\alpha ( 2(n+1)^{r}-n^{r})}
  +\Big(1+\frac{1}{\alpha r(n+2)^{r}}\Big)e^{-\alpha (n+2)^{r}}  \Big),
\end{equation}
where  $O(1)$ is a quantity uniformly bounded in all analyzed parameters.
Formula  (\ref{tel}) contains more exact estimate of remainder in asymptotic decomposition  of the quantity ${\cal E}_{n}(C^{\alpha,r}_{\beta,p})_{C}$ comparing with the estimate  (\ref{step_prepr}).

For $r>1$ and for arbitrary values of  $1\leq p\leq\infty$ the asymptotic equalities for the quantities ${\cal E}_{n}(C^{\alpha,r}_{\beta,p})_{C}$, $\alpha>0$, $\beta\in\mathbb{R}$,
are found in \cite{Serdyuk2005} and have the form
\begin{equation}\label{ser}
  {\cal E}_{n}(C^{\alpha,r}_{\beta,p})_{C}=e^{-\alpha n^{r}}\Big(\frac{\|\cos t\|_{p'}}{\pi}+O(1)\Big(1+\frac{1}{\alpha r n^{r-1}}\Big)e^{-\alpha n^{r-1}}\Big),
\end{equation}
where $O(1)$ is a quantity uniformly bounded in all analyzed parameters. For ${p=\infty}$ the formula (\ref{ser}) follows from (\ref{step_prepr}) and
(\ref{tel}).

 Concerning the case ${0<r<1}$,  except the presented above case  $p=2$, asymptotic equalities for quantities  ${\cal E}_{n}(C^{\alpha,r}_{\beta,p})_{C}$,  $\alpha>0$,  $\beta\in\mathbb{R}$, were known only for $p=\infty$ due to the work of  Stepanets  \cite{Step 1984}, who showed that
\begin{equation}\label{step}
  {\cal E}_{n}(C^{\alpha,r}_{\beta,\infty})_{C}=\frac{4}{\pi^{2}}e^{-\alpha n^{r}}\ln n^{1-r}+O(1)e^{-\alpha n^{r}},
\end{equation}
where $O(1)$ is a quantity uniformly bounded in $n$ and $\beta$.

     In case of  ${0<r<1}$ and  $1\leq p<\infty$ the following order estimates for quantities ${\cal E}_{n}(C^{\alpha,r}_{\beta,p})_{C}$,  $\alpha>0$,  $\beta\in\mathbb{R}$, hold  (see, e.g., \cite{Temlyakov1990Vekya}, \cite{S_S})
\begin{equation}\label{order}
{\cal E}_{n}(C^{\alpha,r}_{\beta,p})_{C}\asymp e^{-\alpha n^{r}}n^{\frac{1-r}{p}}.
\end{equation}

We remark that  for ${0<r<1}$ and $1\leq p<\infty$  Fourier sums provide the order of best approximations of classes $C^{\alpha,r}_{\beta,p}$,  $\alpha>0$,  $\beta\in\mathbb{R}$, in uniform metric, i.e. (see, e.g., [\cite{S_S}, \cite{S_S2})
$$
{\cal E}_{n}(C^{\alpha,r}_{\beta,p})_{C}\asymp{ E}_{n}(C^{\alpha,r}_{\beta,p})_{C}\asymp e^{-\alpha n^{r}}n^{\frac{1-r}{p}},
$$
where
$$
{ E}_{n}(C^{\alpha,r}_{\beta,p})_{C}=\sup\limits_{f\in
C^{\alpha,r}_{\beta,p}}\, \inf\limits_{t_{n-1}\in\mathcal{T}_{2n-1}}\|f-t_{n-1}\|_{C},
$$
and $\mathcal{T}_{2n-1}$ is the subspace of all trigonometric polynomials $t_{n-1}$ of degree not higher than ${n-1}$.

Besides, as  follows from Temlyakov's work   \cite{Temlyakov1990Vekya} for $2\leq p<\infty$ quantities of approximations by Fourier sums realize  order of the linear widths  $\lambda_{2n}$ (definition of  $\lambda_{m}$ see, e.g., \cite[Chapter 1, Section 1.2]{Korn}) of  the classes $C^{\alpha,r}_{0,p}$, i.e.
$$
\lambda_{2n}(C^{\alpha,r}_{0,p},C)\asymp{\cal E}_{n}(C^{\alpha,r}_{\beta,p})_{C}.
$$

In this paper we establish asymptotically sharp estimates of the quantities ${\cal E}_{n}(C^{\alpha,r}_{\beta,p})_{C}$,  ${\alpha>0}$,  $\beta\in\mathbb{R}$,
 for any $0< r<1$ i $1\leq p\leq\infty$. In particular, it is proved, that for  $r\in(0,1)$, $\alpha>0$, $\beta\in\mathbb{R}$ and $1< p\leq\infty$ as $n\rightarrow\infty$ the following asymptotic equality takes place
 $$
 {\cal E}_{n}(C^{\alpha,r}_{\beta,p})_{C}=e^{-\alpha n^{r}}n^{\frac{1-r}{p}}\bigg(
\frac{\|\cos t\|_{p'}}{\pi^{1+\frac{1}{p'}}(\alpha r)^{\frac{1}{p}}}\bigg(\int\limits_{0}^{\infty}\frac{dt}{(t^2+1)^{\frac{p'}{2}}}\bigg)^{\frac{1}{p'}}+
$$
\begin{equation}\label{tan}
O(1)\Big(\frac{1}{n^{(1-r)(p'-1)}}+\frac{1}{n^{r}}+\frac{1}{n^{\frac{1-r}{p}}}\Big)\bigg),
  \end{equation}
   where $\frac{1}{p}+\frac{1}{p'}=1$, and  $O(1)$ is a quantity uniformly bounded with respect to $n$ and $\beta$.
 Herewith, in this paper we found the estimates for remainder in (\ref{tan}),  which are expressed via the parameters of the problem $\alpha, r, p$ in the explicit form and that can be used for practical application.

 \vskip 10mm

 {\bf 2. Formulation of main results.}
For arbitrary $\upsilon>0$ and ${1\leq s\leq \infty}$  assume
\begin{equation}\label{norm_j}
  J_{s}(\upsilon):=\big\|\frac{1}{\sqrt{t^{2}+1}} \big\|_{L_{s}[0,\upsilon]},
\end{equation}
where
$$
\|f \|_{L_{s}[a,b]}=
   {\left\{\begin{array}{cc}
\bigg(\int\limits_{a}^{b}|f(t)|^{s}dt
\bigg)^{\frac{1}{s}}, & 1\leq s<\infty, \\
\mathop{\rm{ess}\sup}\limits_{t\in[a,b]}|f(t)|, \ & s=\infty. \
  \end{array} \right.}
$$
Also for  $\alpha>0$, $r\in(0,1)$ and $1\leq p\leq\infty$ we denote by $n_0=n_0(\alpha,r,p)$ the smallest integer $n$ such that
\begin{equation}\label{n_p}
 \frac{1}{\alpha r}\frac{1}{n^{r}}+\frac{\alpha r \chi(p)}{n^{1-r}}\leq{\left\{\begin{array}{cc}
 \frac{1}{14},  & p=1, \\
\frac{1}{(3\pi)^3}\cdot\frac{p-1}{p}, & 1< p<\infty, \\
\frac{1}{(3\pi)^3}, & p=\infty, \
  \end{array} \right.}
\end{equation}
where $\chi(p)=p$ for $1 \leq p<\infty$ and $\chi(p)=1$ for $p=\infty$.

With the notations introduced above, the main result of this paper is formulated in the following statement:

{\bf Theorem 1.} {\it Let  $0<r<1$, $1\leq p\leq\infty$, $\alpha>0$ and $\beta\in\mathbb{R}$. Then for    $n\geq n_0(\alpha,r,p)$ the following estimate is true $$
{\cal E}_{n}(C^{\alpha,r}_{\beta,p})_{C}=
e^{-\alpha n^{r}}n^{\frac{1-r}{p}}\bigg(\frac{\|\cos t\|_{p'}}{\pi^{1+\frac{1}{p'}}(\alpha r)^{\frac{1}{p}}}J_{p'}\Big(\frac{ \pi n^{1-r}}{\alpha r}\Big)+
$$
\begin{equation}\label{theorem1}
+\gamma_{n,p}^{(1)}\Big(\frac{1}{(\alpha r)^{1+\frac{1}{p}}}J_{p'}\Big(\frac{ \pi n^{1-r}}{\alpha r}\Big)\frac{1}{n^{r}}+\frac{1}{n^{\frac{1-r}{p}}}\Big)\bigg),
\end{equation}
 where $\frac{1}{p}+\frac{1}{p'}=1$, and the quantity ${\gamma_{n,p}^{(1)}=\gamma_{n,p}^{(1)}(\alpha,r,\beta)}$ is such that ${|\gamma_{n,p}^{(1)}|\leq(14\pi)^{2}}$}.

Now we present some corollaries of theorem 1.

For  $1<p<\infty$ theorem 1 yields the following statement:

{\bf Theorem 2.} {\it Let $0<r<1$, $1\leq p<\infty$, $\alpha>0$ and $\beta\in\mathbb{R}$. Then for  $1< p<\infty$ and $n\geq n_0(\alpha,r,p)$ the following estimate is true
$$
 {\cal E}_{n}(C^{\alpha,r}_{\beta,p})_{C}=e^{-\alpha n^{r}}n^{\frac{1-r}{p}}\bigg(
\frac{\|\cos t\|_{p'}}{\pi^{1+\frac{1}{p'}}(\alpha r)^{\frac{1}{p}}}\bigg(\int\limits_{0}^{\infty}\frac{dt}{(t^2+1)^{\frac{p'}{2}}}\bigg)^{\frac{1}{p'}}+
$$
\begin{equation}\label{remark}
+
\gamma^{(2)}_{n,p}\Big(\frac{1}{p'-1}\frac{(\alpha r)^{\frac{p'-1}{p}}}{n^{(1-r)(p'-1)}}+\frac{p^{\frac{1}{p'}}}{(\alpha r)^{1+\frac{1}{p}}}\frac{1}{n^{r}}+\frac{1}{n^{\frac{1-r}{p}}}\Big)\bigg),
\end{equation}
and for $p=1$ and $n\geq n_0(\alpha,r,1)$ the estimate is true
\begin{equation}\label{conseq3}
{\cal E}_{n}(C^{\alpha,r}_{\beta,1})_{C}=
e^{-\alpha n^{r}}n^{1-r}\bigg(\frac{1}{\pi\alpha r}+\gamma_{n,1}^{(2)}\Big(\frac{1}{(\alpha r)^{2}}\frac{1}{n^{r}}+\frac{1}{n^{1-r}}\Big)\bigg),
\end{equation}
 where $\frac{1}{p}+\frac{1}{p'}=1$, and  the quantity ${\gamma_{n,p}^{(2)}=\gamma_{n,p}^{(2)}(\alpha,r,\beta)}$ is such that ${|\gamma_{n,p}^{(2)}|\leq(14\pi)^{2}}$.}

{\bf Proof of the theorem 2.}
According to theorem  1 the following estimate is true for all $1<p<\infty$, $0<r<1$, $\alpha>0$, $\beta\in\mathbb{R}$  and $n\geq n_0(\alpha,r,p)$
$$
{\cal E}_{n}(C^{\alpha,r}_{\beta,p})_{C}=
e^{-\alpha n^{r}}n^{\frac{1-r}{p}}\bigg(\frac{\|\cos t\|_{p'}}{(\alpha r)^{\frac{1}{p}}\pi^{1+\frac{1}{p'}}}\bigg(\int\limits_{0}^{\frac{ \pi n^{1-r}}{\alpha r}}\frac{dt}{(t^2+1)^{\frac{p'}{2}}}\bigg)^{\frac{1}{p'}}+
$$
\begin{equation}\label{theor1}
+\gamma_{n,p}^{(1)}\Big(\frac{1}{(\alpha r)^{1+\frac{1}{p}}}\bigg(\int\limits_{0}^{\frac{ \pi n^{1-r}}{\alpha r}}\frac{dt}{(t^2+1)^{\frac{p'}{2}}}\bigg)^{\frac{1}{p'}}\frac{1}{n^{r}}+\frac{1}{n^{\frac{1-r}{p}}}\Big)\bigg),
\end{equation}
 where $\frac{1}{p}+\frac{1}{p'}=1$, and  the quantity ${\gamma_{n,p}^{(1)}=\gamma_{n,p}^{(1)}(\alpha,r,\beta)}$ is such that  ${|\gamma_{n,p}^{(1)}|\leq(14\pi)^{2}}$.

By applying the Lagrange theorem,  for $n\geq n_0(\alpha,r,p)$ we obtain
 $$
\bigg(\int\limits_{0}^{\infty}\frac{dt}{(t^2+1)^{\frac{p'}{2}}}\bigg)^{\frac{1}{p'}}-
\bigg(\int\limits_{0}^{\frac{ \pi n^{1-r}}{\alpha r}}\frac{dt}{(t^2+1)^{\frac{p'}{2}}}\bigg)^{\frac{1}{p'}}\leq
$$
$$
\leq\frac{1}{p'}
\bigg(\int\limits_{0}^{\frac{ \pi n^{1-r}}{\alpha r}}\frac{dt}{(t^2+1)^{\frac{p'}{2}}}\bigg)^{\frac{1}{p'}-1}
\int\limits_{\frac{\pi n^{1-r}}{\alpha r}}^{\infty}\frac{dt}{(t^2+1)^{\frac{p'}{2}}}
\leq\frac{1}{p'}
\bigg(\int\limits_{0}^{\frac{ \pi n^{1-r}}{\alpha r}}\frac{dt}{(t+1)^{p'}}\bigg)^{-\frac{1}{p}}\int\limits_{\frac{\pi n^{1-r}}{\alpha r}}^{\infty}\frac{dt}{t^{p'}}=
$$
$$
 =\frac{1}{p'}\frac{1}{(p'-1)} \Big(1-\Big(\frac{\pi n^{1-r}}{\alpha r}+1\Big)^{1-p'} \Big)^{-\frac{1}{p}}
 \frac{1}{(p'-1)^{\frac{1}{p}}} \Big(\frac{\alpha r}{\pi n^{1-r}}\Big)^{p'-1}
\leq
$$
$$
 \leq\frac{1}{p'}\frac{1}{(p'-1)} \Big(1-\Big(27\pi^{4}\frac{p^2}{p-1}+1\Big)^{1-p'} \Big)^{-\frac{1}{p}}
 \frac{1}{(p'-1)^{\frac{1}{p}}} \Big(\frac{\alpha r}{\pi n^{1-r}}\Big)^{p'-1}
=
$$
%$$
%=\frac{1}{p'}\frac{1}{(p'-1)^{\frac{1}{p'}}}\Big(\frac{\alpha r}{\pi n^{1-r}}\Big)^{p'-1}
% \Big(1-\Big(27\pi^{4}\frac{p^2}{p-1}+1\Big)^{1-p'} \Big)^{-\frac{1}{p}}=
 %$$
\begin{equation}\label{equation1}
=\frac{1}{p'-1}\Big(\frac{\alpha r}{\pi n^{1-r}}\Big)^{p'-1}
 \frac{(p-1)^{\frac{p-1}{p}}}{p}
\Big(1-\Big(27\pi^{4}\frac{p^2}{p-1}+1\Big)^{\frac{1}{1-p}} \Big)^{-\frac{1}{p}}.
\end{equation}
It can be shown that
$$
 \frac{(p-1)^{\frac{p-1}{p}}}{p}
\Big(1-\Big(27\pi^{4}\frac{p^2}{p-1}+1\Big)^{\frac{1}{1-p}} \Big)^{-\frac{1}{p}}<2.
$$
As follows from (\ref{equation1})
 $$
 \bigg(\int\limits_{0}^{\frac{ \pi n^{1-r}}{\alpha r}}\frac{dt}{(t^2+1)^{\frac{p'}{2}}}\bigg)^{\frac{1}{p'}}=
 $$
\begin{equation}\label{th2}
  =\bigg(\int\limits_{0}^{\infty}\frac{dt}{(t^2+1)^{\frac{p'}{2}}}\bigg)^{\frac{1}{p'}}+
  \frac{\Theta_{\alpha,r,p,n}^{(1)}}{p'-1}\Big(\frac{\alpha r}{\pi n^{1-r}}\Big)^{p'-1}, \ \ |\Theta_{\alpha,r,p,n}^{(1)}|<2.
 \end{equation}
 From relations
 \begin{equation}\label{th21}
\bigg(\int\limits_{0}^{\frac{ \pi n^{1-r}}{\alpha r}}\frac{dt}{(t^2+1)^{\frac{p'}{2}}}\bigg)^{\frac{1}{p'}}\leq
\bigg(\int\limits_{0}^{\infty}\frac{dt}{(t^2+1)^{\frac{p'}{2}}}\bigg)^{\frac{1}{p'}} <
 \bigg(1+\int\limits_{1}^{\infty}\frac{dt}{t^{p'}}\bigg)^{\frac{1}{p'}}<
p^{\frac{1}{p'}}
 \end{equation}
and formulas (\ref{theor1}) and (\ref{th2})  we obtain (\ref{remark}).

Formula (\ref{conseq3}) can be obtained from the equality (\ref{theorem1}) as consequence of substitution $p=1$ and elementary transformations. Theorem 2 is proved.

The following statement follows from the theorem 2 in the case  $p=2$.

{\bf Corollary 1.} {\it Let  $0< r<1$,  $\alpha>0$ and  $\beta\in\mathbb{R}$. Then for $n\geq n_0(\alpha,r,2)$ the following estimate is true
\begin{equation}\label{consequence2}
{\cal E}_{n}(C^{\alpha,r}_{\beta,2})_{C}
 =\frac{e^{-\alpha n^{r}}}{\sqrt{2\pi\alpha r}}n^{\frac{1-r}{2}}\bigg(1+
\gamma^{(1)}_{n}\Big(\frac{1}{\alpha r}\frac{1}{n^{r}}+\frac{\sqrt{\alpha r}}{n^{\frac{1-r}{2}}}\Big)\bigg),
\end{equation}
where the quantity ${\gamma_{n}^{(1)}=\gamma_{n}^{(1)}(\alpha,r,\beta)}$ is such that  ${|\gamma_{n}^{(1)}|\leq392\pi^{\frac{5}{2}}}$.}

{\bf Proof of the corollary 1.} Indeed, setting $p=p'=2$ in the equality (\ref{remark}), we obtain for  $n\geq n_0(\alpha,r,2)$
$$
 {\cal E}_{n}(C^{\alpha,r}_{\beta,2})_{C}=e^{-\alpha n^{r}}n^{\frac{1-r}{2}}\bigg(
\frac{\|\cos t\|_{2}}{\pi^{\frac{3}{2}}(\alpha r)^{\frac{1}{2}}}\bigg(\int\limits_{0}^{\infty}\frac{dt}{t^2+1}\bigg)^{\frac{1}{2}}+
$$
$$
+\gamma^{(2)}_{n,2}\Big(\frac{\sqrt{\alpha r}}{n^{1-r}} +\frac{\sqrt{2}}{(\alpha r)^{\frac{3}{2}}}\frac{1}{n^{r}}+\frac{1}{n^{\frac{1-r}{2}}}\Big)\bigg)=
$$
\begin{equation}\label{triag}
=\frac{e^{-\alpha n^{r}}}{\sqrt{2\pi\alpha r}}n^{\frac{1-r}{2}}\bigg(
1+
\gamma^{(2)}_{n,2}\sqrt{2\pi}\Big(\frac{\alpha r}{n^{1-r}} +\frac{\sqrt{2}}{\alpha r}\frac{1}{n^{r}}+\frac{\sqrt{\alpha r}}{n^{\frac{1-r}{2}}}\Big)\bigg).
\end{equation}

According to (\ref{n_p}) for $n\geq n_0(\alpha,r,2)$
$$
\frac{\sqrt{\alpha r}}{n^{\frac{1-r}{2}}}\leq\frac{1}{2(3\pi)^{\frac{3}{2}}},
$$
therefore
\begin{equation}\label{sta}
\Big(\frac{\alpha r}{n^{1-r}} +\frac{\sqrt{2}}{\alpha r}\frac{1}{n^{r}}+\frac{\sqrt{\alpha r}}{n^{\frac{1-r}{2}}}\Big)\leq
\sqrt{2}\Big(\frac{1}{\alpha r}\frac{1}{n^{r}}+\frac{\sqrt{\alpha r}}{n^{\frac{1-r}{2}}}\Big).
\end{equation}
From (\ref{triag}) and (\ref{sta}) we have (\ref{consequence2}).
Corollary 1 is proved.

However, it is possible to obtain more exact estimate than (\ref{consequence2}) on the basis of equality
 (\ref{serd2011}). Namely, for  $\alpha>0$, $r\in(0,1)$, $\beta\in\mathbb{R}$ and $n\geq n_0(\alpha,r,2)$ the following estimate is true
\begin{equation}\label{sta1}
{\cal E}_{n}(C^{\alpha,r}_{\beta,2})_{C}=\frac{e^{-\alpha n^{r}}}{\sqrt{2\pi\alpha r}}n^{\frac{1-r}{2}}\Big(
1+\gamma_{n}^{(2)}\Big(\frac{1}{2\alpha r}\frac{1}{n^{r}}+\frac{\alpha r}{n^{1-r}}\Big)\Big),
\end{equation}
where  the quantity  ${\gamma_{n}^{(2)}=\gamma_{n}^{(2)}(\alpha,r)}$ is such that $|\gamma_{n}^{(2)}|\leq\sqrt{\frac{54\pi^{3}}{54\pi^{3}-1}} $.
In order to prove (\ref{sta1}) we use the following estimate, which will be useful in what follows.

{\it Let  $\gamma>0$, $r>0$, $m\geq1$ and $\delta\in \mathbb{R}$. Then  for $m\geq\big( \frac{14|\delta+1-r|}{\gamma r}\big)^{\frac{1}{r}}$ the estimate takes place \begin{equation}\label{0lemma_4}
\int\limits_{m}^{\infty}e^{-\gamma t^{r}}t^\delta dt=
\frac{e^{-\gamma m^{r}}}{\gamma r}m^{\delta +1-r}\Big(1+\Theta^{r,\delta}_{\gamma,m}\frac{|\delta+1-r|}{\gamma r}\frac{1}{m^{r}}\Big), \ \ |\Theta^{r,\delta}_{\gamma,m}|\leq\frac{14}{13}.
\end{equation}
}

Indeed,   integrating by parts, we obtain
\begin{equation}\label{0c4}
 \int\limits_{m}^{\infty}e^{-\gamma t^{r}}t^{\delta} dt=\frac{e^{-\gamma m^{r}}}{\gamma r}m^{\delta+1-r}+
\frac{\delta+1-r}{\gamma r}\int\limits_{m}^{\infty}e^{-\gamma t^{r}}t^{-r+\delta} dt.
\end{equation}
Since
\begin{equation}\label{1c4}
\int\limits_{m}^{\infty}e^{-\gamma t^{r}}t^{-r+\delta} dt=\frac{\overline{\Theta}^{r,\delta}_{\gamma,m}}{m^r}\int\limits_{m}^{\infty}e^{-\gamma t^{r}}t^{\delta} dt, \ \ \ 0<\overline{\Theta}^{r,\delta}_{\gamma,m}<1,
\end{equation}
by virtue of (\ref{0c4}) for $m\geq\big( \frac{14|\delta+1-r|}{\gamma r}\big)^{\frac{1}{r}}$ we have
$$
 \int\limits_{m}^{\infty}e^{-\gamma t^{r}}t^\delta dt\leq\frac{e^{-\gamma m^{r}}}{\gamma r}m^{\delta +1-r}+\frac{1}{14}\int\limits_{m}^{\infty}e^{-\gamma t^{r}}t^{\delta} dt,
$$
whence
\begin{equation}\label{0c5}
  \int\limits_{m}^{\infty}e^{-\gamma t^{r}}t^\delta dt\leq\frac{14e^{-\gamma m^{r}}}{13\gamma r}m^{\delta +1-r}.
\end{equation}
The estimate (\ref{0lemma_4}) follows from (\ref{0c4})--(\ref{0c5}).

From the equality (\ref{serd2011}) and relation
 \begin{equation}\label{xi}
\int\limits_{n}^{\infty}\xi(u)du<\sum\limits_{j=n}^{\infty}\xi(j)<\int\limits_{n}^{\infty}\xi(u)du+\xi(n),
\end{equation}
which takes place for any positive and decreasing function $\xi(u)$, $u\geq1$, such that $\int\limits_{n}^{\infty}\xi(u)du<\infty$, we get
  \begin{equation}\label{equation2}
  {\cal E}_{n}(C^{\alpha,r}_{\beta,2})_{C}
=\frac{1}{\sqrt{\pi}}\Big(\int\limits_{n}^{\infty}e^{-2\alpha t^{r}}dt+\Theta_{\alpha,r,n}^{(1)}e^{-2\alpha n^{r}}\Big)^{\frac{1}{2}},  \ \ |\Theta_{\alpha,r,n}^{(1)}|<1.
 \end{equation}

In order to estimate the integral  $\int\limits_{n}^{\infty}e^{-2\alpha t^{r}}dt$ it suffices to use the equality (\ref{0lemma_4}) for $\gamma=2\alpha$, $\delta=0$, $m=n$ and $r\in(0,1)$. Then, taking into account that  ${n_0(\alpha,r,2)>\big( \frac{7(1-r)}{ \alpha r}\big)^{\frac{1}{r}}}$, for $n\geq n_0(\alpha,r,2)$
from (\ref{0lemma_4}) and  (\ref{equation2})  we get
$$
{\cal E}_{n}(C^{\alpha,r}_{\beta,2})_{C}=\frac{1}{\sqrt{\pi}}
\Big(\frac{e^{-2\alpha n^{r}}}{2\alpha r}n^{1-r}\Big(1+\Theta^{r,0}_{2\alpha,n}\frac{(1-r)}{2\alpha r}\frac{1}{n^{r}}\Big)
+\Theta_{\alpha,r,n}^{(1)}e^{-2\alpha n^{r}}\Big)^{\frac{1}{2}}=
$$
\begin{equation}\label{3c5}
=\frac{e^{-\alpha n^{r}}}{\sqrt{2\pi\alpha r}}n^{\frac{1-r}{2}}
\Big(1+\Theta_{\alpha,r,n}^{(2)}\Big(\frac{1}{2\alpha r}\frac{1}{n^{r}}+\frac{\alpha r}{n^{1-r}}\Big)\Big)^{\frac{1}{2}}, \ \ |\Theta_{\alpha,r,n}^{(2)}|\leq2.
\end{equation}
Since for $n>n_0(\alpha,r,2)$
$$
\bigg|\Big(1+\Theta_{\alpha,r,n}^{(2)}\Big(\frac{1}{2\alpha r}\frac{1}{n^{r}}+\frac{\alpha r}{n^{1-r}}\Big)\Big)^{\frac{1}{2}}-1^{\frac{1}{2}}\bigg| \leq
$$
$$
\leq \frac{1}{\sqrt{1-\Big(\frac{1}{\alpha r}\frac{1}{n^{r}}+\frac{2\alpha r}{n^{1-r}}\Big)}}\Big(\frac{1}{2\alpha r}\frac{1}{n^{r}}+\frac{\alpha r}{n^{1-r}}\Big)\leq\sqrt{\frac{54\pi^{3}}{54\pi^{3}-1}} \Big(\frac{1}{2\alpha r}\frac{1}{n^{r}}+\frac{\alpha r}{n^{1-r}}\Big),
$$
then (\ref{sta1}) follows from (\ref{3c5}).

In the case of  $p=\infty$ theorem 1 allows to clarify the asymptotic equality (\ref{step}).

We set $n_1=n_1(\alpha,r)$ be the smallest  number  $n$ such that
\begin{equation}\label{n_2}
 \frac{1}{\alpha r}\frac{1}{n^{r}}\Big(1+\ln\Big(\frac{\pi  n^{1-r}}{\alpha r}\Big)\Big)+\frac{\alpha r }{n^{1-r}}\leq
\frac{1}{(3\pi)^3}.
\end{equation}

The following assertion takes place.

{\bf Theorem 3.} {\it Let  $0< r<1$,  $\alpha>0$  and  $\beta\in\mathbb{R}$. Then for  $n\geq n_1(\alpha,r)$ the following estimate is true
\begin{equation}\label{consequence}
{\cal E}_{n}(C^{\alpha,r}_{\beta,\infty})_{C}=\frac{4}{\pi^{2}}e^{-\alpha n^{r}}\ln\Big(\frac{\pi n^{1-r}}{\alpha r}\Big)+\gamma_{n,\infty}^{(2)}e^{-\alpha n^{r}},
\end{equation}
where  the quantity ${\gamma_{n,\infty}^{(2)}=\gamma_{n,\infty}^{(2)}(\alpha,r,\beta)}$ is such that  ${|\gamma_{n,\infty}^{(2)}|\leq20\pi^{4}}$.}

{\bf Proof of the theorem 3.} From definitions  (\ref{n_2}) and (\ref{n_p}) it follows that ${n_1(\alpha,r)>n_0(\alpha,r,\infty)}$. So, applying the equality (\ref{theorem1}) for $p=\infty \ (p'=1)$,  we get for  $n\geq n_1(\alpha,r)$
\begin{equation}\label{f111}
 {\cal E}_{n}(C^{\alpha,r}_{\beta,\infty})_{C} =
  e^{-\alpha n^{r}}\bigg(
\frac{4}{\pi^{2}}\int\limits_{0}^{\frac{\pi n^{1-r}}{\alpha r}}\frac{dt}{\sqrt{t^{2}+1}}
+\gamma_{n,\infty}^{(1)}\bigg(\frac{1}{\alpha r}\frac{1}{n^{r}}\int\limits_{0}^{\frac{\pi n^{1-r}}{\alpha r}}\frac{dt}{\sqrt{t^{2}+1}}+1\bigg)\bigg).
 \end{equation}

Since
$$
\int\limits_{0}^{\frac{\pi n^{1-r}}{\alpha r}}\frac{dt}{\sqrt{t^{2}+1}}=
\int\limits_{1}^{\frac{\pi n^{1-r}}{\alpha r}}\frac{dt}{t}+
\bigg(\int\limits_{0}^{\frac{\pi n^{1-r}}{\alpha r}}\frac{dt}{\sqrt{t^{2}+1}}-\int\limits_{1}^{\frac{\pi n^{1-r}}{\alpha r}}\frac{dt}{t}\bigg)=
$$
\begin{equation}\label{f112}
 =\ln\Big(\frac{\pi n^{1-r}}{\alpha r}\Big)+\Theta_{\alpha,r,n}^{(3)}, \ \ \ 0<\Theta_{\alpha,r,n}^{(3)}<1,
\end{equation}
by virtue of (\ref{f111}) and (\ref{f112}) for  $n\geq n_1(\alpha,r)$
$$
 {\cal E}_{n}(C^{\alpha,r}_{\beta,\infty})_{C} =
  e^{-\alpha n^{r}}\bigg(
\frac{4}{\pi^{2}}\ln\Big(\frac{\pi n^{1-r}}{\alpha r}\Big)+\frac{4}{\pi^{2}}\Theta_{\alpha,r,n}^{(3)}+
$$
\begin{equation}\label{1f111}
+\gamma_{n,\infty}^{(1)}\bigg(\frac{1}{\alpha r}\frac{1}{n^{r}}\ln\Big(\frac{\pi n^{1-r}}{\alpha r}\Big)+\frac{\Theta_{\alpha,r,n}^{(3)}}{\alpha rn^{r}}+1\bigg)\bigg).
 \end{equation}
The results of our calculations show that for $n\geq n_1(\alpha,r)$
\begin{equation}\label{2f111}
 \frac{4}{\pi^{2}}\Theta_{\alpha,r,n}^{(3)}+
|\gamma_{n,\infty}^{(1)}|\bigg(\frac{1}{\alpha r}\frac{1}{n^{r}}\ln\Big(\frac{\pi n^{1-r}}{\alpha r}\Big)+\frac{\Theta_{\alpha,r,n}^{(3)}}{\alpha rn^{r}}+1\bigg)\leq 20\pi^{4},
 \end{equation}
and therefore, in view of  (\ref{1f111}) and (\ref{2f111}) we obtain (\ref{consequence}).
Theorem 3 is proved.

The asymptotic equality (\ref{step}), which was established by Stepanets, follows from the relation (\ref{consequence}).

{\bf 3.  Proof of the theorem  1.}
According to  (\ref{conv}) and (\ref{sum}) we have
\begin{equation}\label{f1}
{\cal E}_{n}(C^{\alpha,r}_{\beta,p})_{C}=
\frac{1}{\pi}\sup\limits_{\varphi\in B_{p}^{0}}\bigg\|\int\limits_{-\pi}^{\pi}P_{\alpha,r,\beta}^{(n)}(x-t)\varphi(t)dt\bigg\|_{C}, \ \
  \ 1\leq p\leq\infty,
\end{equation}
 where
 \begin{equation}\label{f2}
P_{\alpha,r,\beta}^{(n)}(t):=
\sum\limits_{k=n}^{\infty}e^{-\alpha k^{r}}\cos\Big(kt-\frac{\beta\pi}{2}\Big),  \ 0<r<1, \ \alpha>0, \ \beta\in\mathbb{R}.
\end{equation}
Taking into account the invariance of the sets  $B_{p}^{0}, \ 1\leq p\leq\infty$, under shifts of the argument,
 from (\ref{f1}) we conclude that
\begin{equation}\label{f3}
{\cal E}_{n}(C^{\alpha,r}_{\beta,p})_{C}=
\frac{1}{\pi}\sup\limits_{\varphi\in B_{p}^{0}}\int\limits_{-\pi}^{\pi}P_{\alpha,r,\beta}^{(n)}(t)\varphi(t)dt.
\end{equation}
On the basis of  the duality relation  (see, e.g., \cite[Chapter 1, Section 1.4]{Korn})
\begin{equation}\label{f4}
\sup\limits_{\varphi\in B_{p}^{0}}\int\limits_{-\pi}^{\pi}P_{\alpha,r,\beta}^{(n)}(t)\varphi(t)dt=
\inf\limits_{\lambda\in\mathbb{R}}\|P_{\alpha,r,\beta}^{(n)}(t)-\lambda\|_{p'}, \ \frac{1}{p}+\frac{1}{p'}=1.
\end{equation}

In order to find the estimate for the quantity  $\inf\limits_{\lambda\in\mathbb{R}}\|P_{\alpha,r,\beta}^{(n)}(t)-\lambda\|_{p'}$ we use the following assertion, proof of which will be presented later.

{\bf Lemma $1$.} {\it Let $1\leq s\leq\infty$, $2\pi$--periodic functions $g(t)$ and $h(t)$
have finite derivatives and satisfy the conditions:
\begin{equation}\label{r}
r(t):=\sqrt{g^{2}(t)+h^{2}(t)}\neq0,
\end{equation}
\begin{equation}\label{M}
M:=\sup\limits_{t\in\mathbb{R}}
\frac{\sqrt{(g'(t))^{2}+(h'(t))^{2}}}{\sqrt{g^{2}(t)+h^{2}(t)}}<\infty.
\end{equation}
Then for the function
\begin{equation}\label{snow}
 \phi(t)=g(t)\cos(nt+\gamma)+h(t)\sin(nt+\gamma), \ \ \gamma\in\mathbb{R}, \ \ n\in\mathbb{N},
\end{equation}
for all numbers $ n\geq{\bigg\{\begin{array}{cc}
4\pi s M, & 1\leq s<\infty, \\
1, \ \ \ \ \ \  \ & s=\infty, \
  \end{array} }$
the following estimates take place
\begin{equation}\label{a1}
 \|\phi\|_{s}=\|r\|_{s}\Big(\frac{\|\cos t\|_{s}}{(2\pi)^{\frac{1}{s}}}+\delta_{s,n}^{(1)}\frac{M}{n}\Big),
\end{equation}
\begin{equation}\label{statement2_1}
\inf\limits_{\lambda\in\mathbb{R}}\|\phi(t)-\lambda\|_{s}=\|r\|_{s}\Big(\frac{\|\cos t\|_{s}}{(2\pi)^{\frac{1}{s}}}+\delta_{s,n}^{(2)}\frac{M}{n}\Big),
\end{equation}
\begin{equation}\label{aa2}
 \sup\limits_{h\in\mathbb{R}}\frac{1}{2}\|\phi\big(t+h\big)-\phi(t)\|_{s}= \|r\|_{s}\Big(\frac{\|\cos t\|_{s}}{(2\pi)^{\frac{1}{s}}}+\delta_{s,n}^{(3)}\frac{M}{n}\Big),
 \end{equation}
where}
 \begin{equation}\label{delta}
|\delta_{s,n}^{(i)}|<14\pi, \ i=\overline{1,3}.
 \end{equation}

We represent the function $P_{\alpha,r,\beta}^{(n)}(t)$, which is defined by formula  (\ref{f2}), in the form
\begin{equation}\label{pp}
  P_{\alpha,r,\beta}^{(n)}(t)=
g_{\alpha,r,n}(t)\cos\Big(nt-\frac{\beta\pi}{2}\Big)+h_{\alpha,r,n}(t)\sin\Big(nt-\frac{\beta\pi}{2}\Big),
\end{equation}
where
\begin{equation}\label{g}
  g_{\alpha,r,n}(t):=
\sum\limits_{k=0}^{\infty}e^{-\alpha (k+n)^{r}}\cos kt,
\end{equation}
\begin{equation}\label{h}
 h_{\alpha,r,n}(t):=
-\sum\limits_{k=0}^{\infty}e^{-\alpha (k+n)^{r}}\sin kt.
\end{equation}

Let us show, that for  functions $g_{\alpha,r,n}$ and  $h_{\alpha,r,n}$ the following conditions are satisfied
\begin{equation}\label{cond_lemma}
 \sqrt{g_{\alpha,r,n}^{2}(t)+h_{\alpha,r,n}^{2}(t)}\neq0
\end{equation}
and
\begin{equation}\label{mm1}
 M_{n}=M_{n}(\alpha;r):=\sup\limits_{t\in\mathbb{R}}
 \frac{\sqrt{(g'_{\alpha,r,n}(t))^{2}+(h'_{\alpha,r,n}(t))^{2}}}
 {\sqrt{g_{\alpha,r,n}^{2}(t)+h_{\alpha,r,n}^{2}(t)}}<\infty.
\end{equation}

Since, for arbitrary $\alpha>0, \ 0<r<1$ the sequence  $\big\{ e^{-\alpha (k+n)^{r}}\big\}_{k=0}^{\infty}$ is convex downwards, then  (see, e.g., \cite[Chapter 10, Section 2]{Bari})
$$\frac{1}{2}e^{-\alpha n^{r}}+
\sum\limits_{k=1}^{\infty}e^{-\alpha (k+n)^{r}}\cos kt\geq0,$$
and
\begin{equation}\label{cond01}
  \sqrt{g_{\alpha,r,n}^{2}(t)+h_{\alpha,r,n}^{2}(t)}\geq\frac{1}{2}e^{-\alpha n^{r}}>0.
\end{equation}
 Further, since
  \begin{equation}\label{gg}
  g'_{\alpha,r,n}(t)=
-\sum\limits_{k=1}^{\infty}ke^{-\alpha (k+n)^{r}}\sin kt,
\end{equation}
\begin{equation}\label{hh}
 h'_{\alpha,r,n}(t)=
-\sum\limits_{k=1}^{\infty}ke^{-\alpha (k+n)^{r}}\cos kt,
\end{equation}
it is clear that
 \begin{equation}\label{cond02}
   \sqrt{(g'_{\alpha,r,n}(t))^{2}+(h'_{\alpha,r,n}(t))^{2}}<\sum\limits_{k=1}^{\infty}k e^{-\alpha (k+n)^{r}}<\infty.
\end{equation}
On the basis of (\ref{cond01}) and (\ref{cond02}), the functions $g_{\alpha,r,n}(t)$ and $h_{\alpha,r,n}(t)$ satisfy the conditions  (\ref{cond_lemma}) and (\ref{mm1}).
%В силу твердження 6.6.1 із
%[\ref{Step monog 1987}, p. 253] оскільки послідовність $e^{-\alpha (k+n)^{r}}$, $k\in\mathbb{N}$ опукла донизу і %$\lim\limits_{k\rightarrow\infty}e^{-\alpha (k+n)^{r}}=0$, то
%$g_{\alpha,r,n}^{2}(t)\geq \frac{1}{2}e^{-\alpha n^{r}}$, а отже $g_{\alpha,r,n}^{2}(t)+h_{\alpha,r,n}^{2}(t)\neq0$. Враховуючи цей факт, а також те, що %при фіксованих~$n$ $(g'_{\alpha,r,n}(t))^{2}+(h'_{\alpha,r,n}(t))^{2}<\infty$, отримуємо
Therefore, setting in   lemma 1  $g(t)=g_{\alpha,r,n}(t)$, $h(t)=h_{\alpha,r,n}(t)$, $s=p'$ and $\gamma=-\frac{\beta\pi}{2}$, we get that for
\begin{equation}\label{n_nom}
  n\geq{\bigg\{\begin{array}{cc}
4\pi p' M_{n}, & 1\leq p'<\infty, \\
1, \ \ \ \ \ \  \ & p'=\infty, \
  \end{array} }
\end{equation}
the estimate takes place
\begin{equation}\label{0rrr}
\inf\limits_{\lambda\in\mathbb{R}}\|P_{\alpha,r,\beta}^{(n)}(t)-\lambda\|_{p'}=\Big\| \sqrt{(g_{\alpha,r,n}(t))^{2}+(h_{\alpha,r,n}(t))^{2}}\Big\|_{p'}\Big(\frac{\|\cos t\|_{p'}}{(2\pi)^{\frac{1}{p'}}}+\delta_{n}^{(1)}\frac{M_{n}}{n}\Big),
\end{equation}
where $\frac{1}{p}+\frac{1}{p'}=1$,  quantity $M_{n}$ is defined by equality  (\ref{mm1}), and  the quantity  ${\delta_{n}^{(1)}=\delta_{n}^{(1)}(\alpha,r,\beta,p)}$ is such that  $|\delta_{n}^{(1)}|<14\pi$.

Setting
\begin{equation}\label{r11}
\mathcal{P}_{\alpha,r,n}(t): =g_{\alpha,r,n}(t)-ih_{\alpha,r,n}(t)=\sum\limits_{k=0}^{\infty}e^{-\alpha (k+n)^{r}}e^{ikt},
 \end{equation}
we have
$$
\sqrt{(g'_{\alpha,r,n}(t))^{2}+(h'_{\alpha,r,n}(t))^{2}}=\Big|\mathcal{P'}_{\alpha,r,n}(t) \Big|
$$
and therefore
\begin{equation}\label{mmm1}
 M_{n}=\sup\limits_{t\in\mathbb{R}}
 \frac{\Big|\mathcal{P'}_{\alpha,r,n}(t) \Big|}
 {\Big|\mathcal{P}_{\alpha,r,n}(t) \Big|}.
 \end{equation}

Then, by virtue of the formulas  (\ref{f3}), (\ref{f4}), (\ref{0rrr}) and (\ref{r11}),   for all numbers $n$, which satisfy the condition  (\ref{n_nom}), the estimate holds
\begin{equation}\label{f4440}
  {\cal E}_{n}(C^{\alpha,r}_{\beta,p})_{C}=\|\mathcal{P}_{\alpha,r,n}(t)\|_{p'}
\Big(\frac{\|\cos t\|_{p'}}{2^{\frac{1}{p'}}\pi^{1+\frac{1}{p'}}}
+\delta_{n}^{(2)}\frac{M_{n}}{n}\Big), \ 1\leq p\leq\infty,
\end{equation}
 where  $M_{n}$ is defined by equality (\ref{mmm1}), and for the quantity  ${\delta_{n}^{(2)}=\delta_{n}^{(2)}(\alpha,r,\beta,p)}$  is such that $|\delta_{n}^{(2)}|<14$.

Since
\begin{equation}\label{sqr_p}
 \Big|\mathcal{P}_{\alpha,r,n}(t) \Big|^{2}=\mathcal{P}_{\alpha,r,n}(t)\widetilde{\mathcal{P}}_{\alpha,r,n}(t),
\end{equation}
where
 $$
 \widetilde{\mathcal{P}}_{\alpha,r,n}(t)=g_{\alpha,r,n}(t)+ih_{\alpha,r,n}(t)=\sum\limits_{k=0}^{\infty}e^{-\alpha (k+n)^{r}}e^{-ikt},
 $$
  by expanding the product
$\mathcal{P}_{\alpha,r,n}\widetilde{\mathcal{P}}_{\alpha,r,n}$ in the Fourier series (see, e.g., \cite[Chapter~1, Section 23]{Bari}), we get
$$
\mathcal{P}_{\alpha,r,n}(t)\widetilde{\mathcal{P}}_{\alpha,r,n}(t)=
\Big(\sum\limits_{k=0}^{\infty}e^{-\alpha (k+n)^{r}}e^{ikt}\Big)
\Big(\sum\limits_{k=-\infty}^{0}e^{-\alpha (-k+n)^{r}}e^{ikt}\Big)=
$$
$$
 =\sum\limits_{k=-\infty}^{\infty}\sum\limits_{j=0}^{\infty}e^{-\alpha (j+n)^{r}}e^{-\alpha (j+|k|+n)^{r}}e^{ikt}=
$$
\begin{equation}\label{f86}
=
\sum\limits_{j=n}^{\infty}e^{-2\alpha j^{r}}+
2\sum\limits_{k=1}^{\infty}\sum\limits_{j=n}^{\infty}e^{-\alpha j^{r}}e^{-\alpha (j+k)^{r}}\cos kt.
\end{equation}

Let convert the sum  $\sum\limits_{j=n}^{\infty}e^{-2\alpha j^{r}}+
2\sum\limits_{k=1}^{\infty}\sum\limits_{j=n}^{\infty}e^{-\alpha j^{r}}e^{-\alpha (j+k)^{r}}\cos kt$ with a help of Poisson summation formula.

{\bf Assertion 1 \cite[Chapter 2, Section 2.8]{Titmarsh}.} {\it
 Let continuous function  $\phi(x)$ be a function of bounded variation
in the interval $(0,\infty)$, $\lim\limits_{x\rightarrow\infty}\phi(x)=0$
and
$$
\int\limits_{0}^{\infty}\phi(t)dt<\infty.
$$
Then the following equality takes place
\begin{equation}\label{st}
\sqrt{a}\Big(\frac{\phi(0)}{2}+\sum\limits_{k=1}^{\infty}\phi(ka)\Big)=
\sqrt{\frac{2\pi}{a}}\Big(\frac{\Phi_{c}(0)}{2}+\sum\limits_{k=1}^{\infty}\Phi_{c}(\frac{2\pi k}{a})\Big), \ \ a>0,
\end{equation}
where $\Phi_{c}(x)$ is the Fourier cosine transform of the function  $\phi(x)$ of the form}
$$
\Phi_{c}(x)=\sqrt{\frac{2}{\pi}}\int\limits_{0}^{\infty}\phi(u)\cos xu du.
$$

Let fix $t\in[-\pi, \pi]$, $\alpha>0$, ${r\in(0,1)}$ and set
$$\phi(x)=2\sum\limits_{j=n}^{\infty}e^{-\alpha j^{r}}e^{-\alpha (j+x)^{r}}\cos xt,\ \ \ \ x\geq0$$
   and $a=1$. One can easily check that all conditions of the assertion 1 are satisfied, and therefore, according to (\ref{st}) we obtain
$$
\sum\limits_{j=n}^{\infty}e^{-2\alpha j^{r}}+
2\sum\limits_{k=1}^{\infty}\sum\limits_{j=n}^{\infty}e^{-\alpha j^{r}}e^{-\alpha (j+k)^{r}}\cos kt=
$$
$$
=2\int\limits_{0}^{\infty}\sum\limits_{j=n}^{\infty}e^{-\alpha j^{r}}e^{-\alpha (j+u)^{r}}\cos ut du+
$$
$$
+4\sum\limits_{k=1}^{\infty}\int\limits_{0}^{\infty}\sum\limits_{j=n}^{\infty}e^{-\alpha j^{r}}e^{-\alpha (j+u)^{r}}\cos ut\cos 2\pi ku du=
$$
\begin{equation}\label{f105}
  =Q_{n}(t)+R_{n}(t),
\end{equation}
where
\begin{equation}\label{rn1}
 Q_{n}(t)=Q_{n}(\alpha; r; t):=2\sum\limits_{j=n}^{\infty}e^{-\alpha j^{r}}\int\limits_{0}^{\infty}e^{-\alpha (j+u)^{r}}\cos ut du,
\end{equation}
$$
R_{n}(t)=R_{n}(\alpha; r; t):=
$$
\begin{equation}\label{0rn1}
:=2\sum\limits_{k=1}^{\infty}\sum\limits_{j=n}^{\infty}e^{-\alpha j^{r}}
\int\limits_{0}^{\infty}e^{-\alpha (j+u)^{r}}\big(\cos((t-2\pi k)u)+\cos((t+2\pi k)u)\big)du.
\end{equation}

Hence, as a consequence of (\ref{sqr_p}), (\ref{f86}) and (\ref{f105})
\begin{equation}\label{b22}
\Big|\mathcal{P}_{\alpha,r,n}(t) \Big|^{2}=Q_{n}(t)+R_{n}(t).
\end{equation}

Denote by  $n_2=n_2(\alpha,r,p)$ the smallest  number  $n$ such that
\begin{equation}\label{n1}
\frac{1}{\alpha r}\frac{1}{n^{r}}+\frac{\alpha r \chi(p)}{n^{1-r}}\leq\frac{1}{14},
\end{equation}
where
$$
\chi(p)={\bigg\{\begin{array}{cc}
p, & 1\leq p<\infty, \\
1,  & p=\infty, \
  \end{array} }
$$
and let us show that for the quantity  $ Q_{n}(t)$ for  $n\geq n_{2}(\alpha,r,p)$ and arbitrary $t\in[-\pi, \pi]$ the following estimate takes place
\begin{equation}\label{f34}
Q_{n}(t)=\frac{ e^{-2\alpha n^{r}}}{t^{2}+(\alpha r n^{r-1})^{2}}\Big(
1+\Theta_{\alpha,r,n}^{(4)}(t)\Big(\frac{1-r}{\alpha r}\frac{1}{n^{r}}+\frac{\alpha r}{n^{1-r}}\Big)\Big), \ \ |\Theta_{\alpha,r,n}^{(4)}(t)|<5.
\end{equation}

Integrating by parts, we find
$$
\int e^{-\alpha (j+u)^{r}}\cos utdu=
$$
$$
=e^{-\alpha (j+u)^{r}}\frac{-\alpha r(j+u)^{r-1}\cos ut+t\sin ut}
{t^{2}+(\alpha r(j+u)^{r-1})^{2}}+\alpha r(1-r)\times
$$
$$
\times\int e^{-\alpha (j+u)^{r}}(j+u)^{r-2}\frac{((\alpha r(j+u)^{r-1})^{2}-t^{2})\cos ut-
2t\alpha r(j+u)^{r-1}\sin ut
}{(t^{2}+(\alpha r(j+u)^{r-1})^{2})^{2}}du.
$$

Hence, we obtain the equality
$$
\int\limits_{0}^{\infty} e^{-\alpha (j+u)^{r}}\cos utdu=
\frac{\alpha rj^{r-1}}{t^{2}+(\alpha rj^{r-1})^{2}}e^{-\alpha j^{r}}+\alpha r(1-r)\times
$$
\begin{equation}\label{c1}
  \times
\int\limits_{0}^{\infty}e^{-\alpha (j+u)^{r}}(j+u)^{r-2}\frac{((\alpha r(j+u)^{r-1})^{2}-t^{2})\cos ut-
2t\alpha r(j+u)^{r-1}\sin ut
}{(t^{2}+(\alpha r(j+u)^{r-1})^{2})^{2}}du.
\end{equation}

It is easy to verify that
$$
\bigg|
\int\limits_{0}^{\infty}e^{-\alpha (j+u)^{r}}(j+u)^{r-2}\frac{((\alpha r(j+u)^{r-1})^{2}-t^{2})\cos ut-
2t\alpha r(j+u)^{r-1}\sin ut
}{((\alpha r(j+u)^{r-1})^{2}+t^{2})^{2}}du\bigg|\leq
$$
$$
\leq\int\limits_{0}^{\infty}e^{-\alpha (j+u)^{r}}(j+u)^{r-2}\bigg(\frac{1
}{t^{2}+(\alpha r(j+u)^{r-1})^{2}}+
\frac{2t\alpha r(j+u)^{r-1}}{(t^{2}+(\alpha r(j+u)^{r-1})^{2})^{2}}\bigg)du\leq
$$
\begin{equation}\label{c2}
 \leq2
\int\limits_{0}^{\infty}e^{-\alpha (j+u)^{r}}\frac{(j+u)^{r-2}
}{t^{2}+(\alpha r(j+u)^{r-1})^{2}}du.
\end{equation}
For fixed $\alpha>0, \ r\in(0,1)$ and ${t\in[-\pi,\pi]}$ the function $\frac{v^{r-2}
}{t^{2}+(\alpha rv^{r-1})^{2}}, \ v\geq1$ decreases. Besides, according to   (\ref{0c5}), for $\delta=0$, $\gamma=\alpha$, $m=j$, $j\geq n_{2}(\alpha,r,p)$ the estimate takes place
$$
\int\limits_{0}^{\infty}e^{-\alpha (j+u)^{r}}\frac{(j+u)^{r-2}
}{t^{2}+(\alpha r(j+u)^{r-1})^{2}}du \leq \frac{j^{r-2}
}{t^{2}+(\alpha rj^{r-1})^{2}}\int\limits_{0}^{\infty}e^{-\alpha (j+u)^{r}}du=
$$
\begin{equation}\label{c3}
=\frac{j^{r-2}
}{t^{2}+(\alpha rj^{r-1})^{2}}\int\limits_{j}^{\infty}e^{-\alpha u^{r}}du\leq\frac{14}{13}\frac{e^{-\alpha j^{r}}
}{\alpha r j(t^{2}+(\alpha rj^{r-1})^{2})}.
\end{equation}

It follows from relations  (\ref{c1})--(\ref{c3}) that for $j\geq n_{2}(\alpha,r,p)$
$$
\int\limits_{0}^{\infty} e^{-\alpha (j+u)^{r}}\cos utdu=
$$
\begin{equation}\label{f100}
=\frac{\alpha rj^{r-1}}{t^{2}+(\alpha rj^{r-1})^{2}}e^{-\alpha j^{r}}\Big(1+\Theta_{\alpha,r,j}^{(5)}(t)\frac{1-r}{\alpha r}\frac{1}{j^{r}}\Big), \ \ |\Theta_{\alpha,r,j}^{(5)}(t)|\leq\frac{28}{13}.
\end{equation}
Therefore, taking into account  (\ref{rn1}),  for $n\geq n_{2}(\alpha,r,p)$ we have
\begin{equation}\label{f101}
  Q_{n}(t)=2\alpha r\sum\limits_{j=n}^{\infty}\frac{e^{-2\alpha j^{r}}j^{r-1}}{t^{2}+(\alpha rj^{r-1})^{2}}\Big(1+\Theta_{\alpha,r,n}^{(6)}(t)\frac{1-r}{\alpha r}\frac{1}{n^{r}}\Big), \ \ |\Theta_{\alpha,r,n}^{(6)}(t)|\leq\frac{28}{13}.
\end{equation}

Further, let us find bilateral estimates for the quantities  $\sum\limits_{j=n}^{\infty}\frac{e^{-2\alpha j^{r}} j^{r-1}}{t^{2}+(\alpha r j^{r-1})^{2}}$ for ${n\geq n_2(\alpha,r,p)}$.
  It can be shown that for fixed   $\alpha>0$, ${r\in(0,1)}$ and $t\in[-\pi, \pi]$ the function  ${\xi(u)=\frac{e^{-2\alpha u^{r}} u^{r-1}}{t^{2}+(\alpha r u^{r-1})^{2}}}$ decreases for $u\geq n_{2}(\alpha,r,p)$. Therefore, on basis of  (\ref{xi})
$$
 2\alpha r\sum\limits_{j=n}^{\infty}e^{-2\alpha j^{r}}\frac{ j^{r-1}}{t^{2}+(\alpha r j^{r-1})^{2}}=
$$
\begin{equation}\label{f21}
=
2\alpha r\int\limits_{n}^{\infty}e^{-2\alpha u^{r}}\frac{ u^{r-1}}{t^{2}+(\alpha r u^{r-1})^{2}}du
+\Theta_{\alpha,r,n}^{(7)}(t)\frac{\alpha r \ e^{-2\alpha n^{r}} n^{r-1}}{t^{2}+(\alpha r n^{r-1})^{2}}, \ \ 0\leq\Theta_{\alpha,r,n}^{(7)}(t)\leq2.
\end{equation}
Integrating by parts, we have
$$
2\alpha r\int\limits_{n}^{\infty}\frac{ e^{-2\alpha u^{r}}u^{r-1}}{t^{2}+(\alpha r u^{r-1})^{2}}du=
$$
\begin{equation}\label{c8}
 =
\frac{ e^{-2\alpha n^{r}}}{t^{2}+(\alpha r n^{r-1})^{2}}+
2(\alpha r)^{2}(1-r)\int\limits_{n}^{\infty}\frac{ e^{-2\alpha u^{r}}u^{2r-3}}{(t^{2}+(\alpha r u^{r-1})^{2})^{2}}du.
\end{equation}

Since
$$
(\alpha r)^{2}\int\limits_{n}^{\infty}\frac{e^{-2\alpha u^{r}} u^{2r-3}}{(t^{2}+(\alpha r u^{r-1})^{2})^{2}}du\leq
\int\limits_{n}^{\infty}\frac{e^{-2\alpha u^{r}} u^{-1}}{t^{2}+(\alpha r u^{r-1})^{2}}du\leq
$$
\begin{equation}\label{rrr}
\leq\frac{1}{n^{r}}\int\limits_{n}^{\infty}\frac{e^{-2\alpha u^{r}} u^{r-1}}{t^{2}+(\alpha r u^{r-1})^{2}}du,
\end{equation}
it follows from  (\ref{c8}) that for $n\geq n_{2}(\alpha,r,p)$ the following inequalities are true
$$
\int\limits_{n}^{\infty}\frac{ e^{-2\alpha u^{r}}u^{r-1}}{t^{2}+(\alpha r u^{r-1})^{2}}du\leq
$$
$$
\leq
\frac{1}{2\alpha r}\frac{ e^{-2\alpha n^{r}}}{t^{2}+(\alpha r n^{r-1})^{2}}+
\frac{1-r}{\alpha r}\frac{1}{n^{r}}\int\limits_{n}^{\infty}\frac{e^{-2\alpha u^{r}} u^{r-1}}{t^{2}+(\alpha r u^{r-1})^{2}}du\leq
$$
$$
\leq
\frac{1}{2\alpha r}\frac{ e^{-2\alpha n^{r}}}{t^{2}+(\alpha r n^{r-1})^{2}}+
\frac{1}{14}\int\limits_{n}^{\infty}\frac{e^{-2\alpha u^{r}} u^{r-1}}{t^{2}+(\alpha r u^{r-1})^{2}}du.
$$
Hence, for $n\geq n_{2}(\alpha,r,p)$
\begin{equation}\label{c9}
 \int\limits_{n}^{\infty}\frac{e^{-2\alpha u^{r}} u^{r-1}}{t^{2}+(\alpha r u^{r-1})^{2}}du\leq\frac{7}{13\alpha r}
\frac{ e^{-2\alpha n^{r}}}{t^{2}+(\alpha r n^{r-1})^{2}}.
\end{equation}
From (\ref{c8})--(\ref{c9}) for $ n\geq n_{2}(\alpha,r,p)$ we arrive at the following estimate
$$
 2\alpha r\int\limits_{n}^{\infty}e^{-2\alpha u^{r}}\frac{ u^{r-1}}{t^{2}+(\alpha r u^{r-1})^{2}}du=
$$
\begin{equation}\label{f24}
=
 \frac{ e^{-2\alpha n^{r}}}{t^{2}+(\alpha r n^{r-1})^{2}}\Big(
1+\Theta_{\alpha,r,n}^{(8)}(t)\frac{1-r}{\alpha r}\frac{1}{n^{r}}\Big), \ \ 0<\Theta_{\alpha,r,n}^{(8)}(t)\leq\frac{14}{13}.
\end{equation}
It follows from formulas  (\ref{f21}) and (\ref{f24}) that
$$
2\alpha r \sum\limits_{j=n}^{\infty}\frac{e^{-2\alpha j^{r}} j^{r-1}}{t^{2}+(\alpha r j^{r-1})^{2}}=
$$
\begin{equation}\label{f25}
=
\frac{ e^{-2\alpha n^{r}}}{t^{2}+(\alpha r n^{r-1})^{2}}\big(
1+\Theta_{\alpha,r,n}^{(9)}(t)\Big(\frac{1-r}{\alpha r}\frac{1}{n^{r}}+\frac{\alpha r}{n^{1-r}}\Big)\Big),
\end{equation}
where $ n\geq n_{2}(\alpha,r,p)$ and $0<\Theta_{\alpha,r,n}^{(9)}(t)\leq2$.

In view of  (\ref{f101}) and (\ref{f25}) for all $n\geq n_{2}(\alpha,r,p)$ we obtain (\ref{f34}). In particular, it follows from formulas   (\ref{n1}) and (\ref{f34}) that
\begin{equation}\label{dod1}
Q_{n}(t)>0, \ t\in[-\pi, \pi], \  n\geq n_{2}(\alpha,r,p).
\end{equation}

Let us find upper estimate for the quantity  $R_{n}(t)$ of the form (\ref{0rn1}).
Denote by  ${\mathfrak M}$ the set of all convex downwards, continuous functions  $\psi(t)>0, \ t\geq 1,$ such that $\lim\limits_{t\rightarrow\infty}\psi(t)=0$.
The following assertion takes place.

{\bf Lemma 2.} {\it Let $\psi\in\mathfrak{M}$. Then }
\begin{equation}\label{f10}
0<\int\limits_{0}^{\infty}\psi(\tau+u)\cos v u du
\leq \frac{\pi}{v^{2}}|\psi'(\tau)|, \ \ v\in\mathbb{R}\setminus\{0\}, \ \ \tau\geq1.
\end{equation}

{\bf Proof of lemma 2.} We use the scheme of the proof of the estimate (2.4.31) from the work \cite[p.~93]{St-Ruk-Ch}.
Let, e.g., consider the case $v>0$. Using the method of integration by parts, we have
\begin{equation}\label{f46}
\int\limits_{0}^{\infty}\psi(\tau+u)\cos v u du
=\frac{-1}{v}\int\limits_{0}^{\infty}\psi'(\tau+u)\sin v u du.
\end{equation}
We set
$$
I(x)=I(\psi;\tau;v;x):=-\int\limits_{x}^{\infty}\psi'(j+u)\sin v u du, \ x\geq0, \ v>0, \ \tau\in\mathbb{N}.
$$

The function $I(x)$, obviously, is continuous for every fixed   $v$, and on every interval
between the consecutive zeros  $u_{m}=\frac{\pi m}{v}$ and $u_{m+1}=\frac{\pi (m+1)}{v}$
of the function $\sin vu$ has one simple zero  $x_{m}$. Existance of zeros  $x_{m}$ of the function $I(x)$ is a consequence of the Leibniz theorem on
alternating series, and uniqueness of zero   $x_{m}$ on the interval  $(u_{m}, u_{m+1})$ follows from the
equality
$$
\mathrm{sign}\, I'(x)=-\mathrm{sign} \sin xv, \ x\in(u_{m}, u_{m+1})\ m\in\mathbb{Z}_{+}.
$$

Let $x_{0}$ be the zero closest from the right to the point  $x=0$. It is obvious that
$$
0\leq x_{0}\leq \frac{\pi}{v}.
$$
Taking into account this fact and also monotone decreasing of the function  $-\psi'(t)$ on the interval $[1,\infty)$, we have
$$
\frac{-1}{v}\int\limits_{0}^{\infty}\psi'(\tau+u)\sin v u du
=\frac{1}{v}\int\limits_{0}^{x_{0}}|\psi'(\tau+u)|\sin v u du\leq
$$
\begin{equation}\label{f31}
  \leq\frac{1}{v}\int\limits_{0}^{\frac{\pi}{v}}|\psi'(\tau+u)| du
\leq \frac{\pi}{v^{2}}|\psi'(\tau)|.
\end{equation}

For $v>0$  inequality (\ref{f10}) follows from the formulas  (\ref{f46}) and (\ref{f31}).
For $v<0$  the proof of inequality  (\ref{f10}) is analogous. Lemma 2 is proved.

Setting in  inequality (\ref{f10})  $v=t\pm2\pi k, \ k\in\mathbb{N}$, and $\tau=j$,  we obtain that for arbitrary $\psi\in\mathfrak{M}$ and $t\in[-\pi,\pi]$
$$
0<\sum\limits_{k=1}^{\infty}\sum\limits_{j=n}^{\infty}\psi(j)\int\limits_{0}^{\infty}\psi(j+u)\big(\cos((t-2\pi k)u)+\cos((t+2\pi k)u)\big)du\leq
$$
$$
\leq\pi
\sum\limits_{k=1}^{\infty}
\Big(\frac{1}{(t-2k\pi)^{2}}+\frac{1}{(t+2k\pi)^{2}}\Big)\sum\limits_{j=n}^{\infty}\psi(j)|\psi'(j)|\leq
$$
$$
\leq\pi
\sum\limits_{k=1}^{\infty}\Big(\frac{1}{(\pi-2k\pi)^{2}}+
\frac{1}{(\pi+2k\pi)^{2}}\Big)\psi(n)\bigg(|\psi'(n)|+\int\limits_{n}^{\infty}|\psi'(u)|du\bigg)=
$$
$$
=\frac{1}{\pi}
\sum\limits_{k=1}^{\infty}\Big(\frac{1}{(2k-1)^{2}}+\frac{1}{(2k+1)^{2}}\Big)
\psi(n)\bigg(|\psi'(n)|+\psi(n)\bigg)=
$$
\begin{equation}\label{lemma_3}
=\Big(\frac{\pi}{4}-\frac{1}{\pi}\Big)
\psi(n)\bigg(|\psi'(n)|+\psi(n)\bigg), \ n\in\mathbb{N}.
\end{equation}

Setting in    (\ref{lemma_3})    $\psi(t)=e^{-\alpha t^{r}}$, ${0<r< 1}$, $\alpha>0$, we get that  for the function $R_{n}(t)$ of the form (\ref{0rn1}) the following estimate takes place
\begin{equation}\label{f15}
0<R_{n}(t)\leq \Big(\frac{\pi}{2}-\frac{2}{\pi}\Big)e^{-2\alpha n^{r}}\big(\frac{\alpha r}{n^{1-r}}+1\big)\leq \Big(\frac{\pi}{2}-\frac{2}{\pi}\Big)\frac{15}{14}e^{-2\alpha n^{r}}<\frac{\pi}{3}e^{-2\alpha n^{r}},
\end{equation}
where $n\geq n_{2}(\alpha,r,p)$.

By virtue of (\ref{b22})
\begin{equation}\label{0p}
\big |\mathcal{P}_{\alpha,r,n}(t) \big|=
 \sqrt{Q_{n}(t)+R_{n}(t)},
\end{equation}
and therefore,  taking into account  (\ref{dod1}) and (\ref{f15}),  we have
\begin{equation}\label{f104}
\big \|\mathcal{P}_{\alpha,r,n} \big\|_{p'}=
 \big\|\sqrt{Q_{n}}\big\|_{L_{p'}[-\pi,\pi]}+\Theta_{\alpha,r,p,n}^{(2)}e^{-\alpha n^{r}}, \ \ \ 1\leq p'\leq\infty,
\end{equation}
where $|\Theta_{\alpha,r,p,n}^{(2)}|<\frac{2\pi^{2}}{3}$ and $ n\geq n_{2}(\alpha,r,p)$.

Let us show, that for $1\leq p'\leq\infty$, $\frac{1}{p}+\frac{1}{p'}=1$, and $n\geq n_2(\alpha,r,p)$ the estimate is true
$$
\big\|\mathcal{P}_{\alpha,r,n} \big\|_{p'}=e^{-\alpha n^{r}}n^{\frac{1-r}{p}}\bigg(\frac{2^{\frac{1}{p'}}}{(\alpha r)^{\frac{1}{p}}}J_{p'}\Big(\frac{ \pi n^{1-r}}{\alpha r}\Big)+
$$
\begin{equation}\label{f109}
 +\Theta_{\alpha,r,p,n}^{(3)}\Big(\frac{1-r}{(\alpha r)^{1+\frac{1}{p}}}J_{p'}\Big(\frac{ \pi n^{1-r}}{\alpha r}\Big)\frac{1}{n^{r}}+\frac{1}{n^{\frac{1-r}{p}}}\Big)\bigg),
\end{equation}
where
 \begin{equation}\label{theta3}
 |\Theta_{\alpha,r,p,n}^{(3)}|\leq {\left\{\begin{array}{cc}
\pi^{2}, & 1\leq p'<\infty, \\
\frac{14}{13}, \ & p'=\infty. \
  \end{array} \right.}
\end{equation}

Since, on the basis of estimate  (\ref{f34})  for  $n\geq n_{2}(\alpha,r,p)$ and $1\leq p'\leq\infty$
$$
\Big|\Big(
1+\Theta_{\alpha,r,n}^{(10)}(t)\Big(\frac{1-r}{\alpha r}\frac{1}{n^{r}}+\frac{\alpha r}{n^{1-r}}\Big)^{\frac{1}{2}}-1\Big|\leq
$$
$$
\leq\frac{5}{2}\frac{1}{\sqrt{1-5\Big(\frac{1-r}{\alpha r}\frac{1}{n^{r}}+\frac{\alpha r}{n^{1-r}}\Big)}}\Big(\frac{1-r}{\alpha r}\frac{1}{n^{r}}+\frac{\alpha r}{n^{1-r}}\Big)\leq \frac{5\sqrt{7}}{3\sqrt{2}}\Big(\frac{1-r}{\alpha r}\frac{1}{n^{r}}+\frac{\alpha r}{n^{1-r}}\Big)
$$
we get
\begin{equation}\label{f3334}
\sqrt{Q_{n}(t)}=\frac{ e^{-\alpha n^{r}}}{\sqrt{t^{2}+(\alpha r n^{r-1})^{2}}}\Big(
1+\Theta_{\alpha,r,n}^{(10)}(t)\Big(\frac{1-r}{\alpha r}\frac{1}{n^{r}}+\frac{\alpha r}{n^{1-r}}\Big)\Big),  \ \ |\Theta_{\alpha,r,n}^{(10)}(t)|\leq \frac{5\sqrt{7}}{3\sqrt{2}}.
\end{equation}
For $1\leq p'<\infty$ from (\ref{f3334}) we have
$$
\Big\|\sqrt{Q_{n}}\Big\|_{L_{p'}[-\pi,\pi]}=e^{-\alpha n^{r}}\bigg(\int\limits_{-\pi}^{\pi}\frac{dt}{(t^{2}+(\alpha r n^{r-1})^{2})^{\frac{p'}{2}}}\bigg)^{\frac{1}{p'}}
\Big(
1+\Theta_{\alpha,r,p,n}^{(4)}\Big(\frac{1-r}{\alpha r}\frac{1}{n^{r}}+\frac{\alpha r}{n^{1-r}}\Big)\Big)=
$$
%$$
%=2^{\frac{1}{p'}}e^{-\alpha n^{r}}\bigg(\int\limits_{0}^{\pi}\frac{dt}{(t^{2}+(\alpha r n^{r-1})^{2})^{\frac{p'}{2}}}\bigg)^{\frac{1}{p'}}
%\Big(
%1+\Theta_{\alpha,r,p,n}^{(4)}\Big(\frac{1-r}{\alpha r}\frac{1}{n^{r}}+\frac{\alpha r}{n^{1-r}}\Big)\Big)=
%$$
\begin{equation}\label{f102}
  =2^{\frac{1}{p'}}e^{-\alpha n^{r}}\Big(\frac{ n^{1-r}}{\alpha r}\Big)^{\frac{1}{p}}J_{p'}\Big(\frac{\pi n^{1-r}}{\alpha r}\Big)
\Big(
1+\Theta_{\alpha,r,p,n}^{(4)}\Big(\frac{1-r}{\alpha r}\frac{1}{n^{r}}+\frac{\alpha r}{n^{1-r}}\Big)\Big),
\end{equation}
where $|\Theta_{\alpha,r,p,n}^{(4)}|\leq \frac{5\sqrt{7}}{3\sqrt{2}}$ and $J_{p'}\Big(\frac{\pi n^{1-r}}{\alpha r}\Big)$ is defined by equality (\ref{norm_j}).

Combining  (\ref{f104}) and (\ref{f102}), we obtain that for  $1\leq p'<\infty$ the following relation takes place
$$
\big\|\mathcal{P}_{\alpha,r,n} \big\|_{p'}=e^{-\alpha n^{r}}n^{\frac{1-r}{p}}\bigg(\frac{2^{\frac{1}{p'}}}{(\alpha r)^{\frac{1}{p}}}J_{p'}\Big(\frac{ \pi n^{1-r}}{\alpha r}\Big)+
$$
\begin{equation}\label{f200}
 +\Theta_{\alpha,r,p,n}^{(4)}\frac{2^{\frac{1}{p'}}}{(\alpha r)^{\frac{1}{p}}}J_{p'}\Big(\frac{ \pi n^{1-r}}{\alpha r}\Big)\Big(\frac{1-r}{\alpha r}\frac{1}{n^{r}}+\frac{\alpha r}{n^{1-r}}\Big)+\frac{\Theta_{\alpha,r,p,n}^{(2)}}{n^{\frac{1-r}{p}}}\bigg).
\end{equation}

However, for all $n>n_2(\alpha,r,p)$
\begin{equation}\label{f201}
 \frac{2^{\frac{1}{p'}}}{(\alpha r)^{\frac{1}{p}}}J_{p'}\Big(\frac{ \pi n^{1-r}}{\alpha r}\Big)
 \frac{\alpha r}{n^{1-r}}<\frac{1}{n^{\frac{1-r}{p}}}, \ 1\leq p'<\infty.
\end{equation}

Indeed, taking into account  (\ref{norm_j}) and (\ref{n1}), for all  $ 1< p'<\infty$ and ${n\geq n_2(\alpha,r,p)}$ we find
$$
 \frac{2^{\frac{1}{p'}}}{(\alpha r)^{\frac{1}{p}}}J_{p'}\Big(\frac{ \pi n^{1-r}}{\alpha r}\Big)
 \frac{\alpha r}{n^{1-r}}n^{\frac{1-r}{p}}=
 \Big(\frac{2\alpha r}{n^{1-r}}\Big)^{\frac{1}{p'}}J_{p'}\Big(\frac{ \pi n^{1-r}}{\alpha r}\Big)
  <
$$
$$
   <\Big(\frac{2\alpha r}{n^{1-r}}\Big)^{\frac{1}{p'}}\bigg(\int\limits_{0}^{\infty}\frac{dt}{(t^{2}+1)^{\frac{p'}{2}}}\bigg)^{\frac{1}{p'}}
<\Big(\frac{2\alpha r}{n^{1-r}}\Big)^{\frac{1}{p'}}\bigg(1+\int\limits_{1}^{\infty}\frac{dt}{t^{p'}}\bigg)^{\frac{1}{p'}}=
$$
\begin{equation}\label{in1}
  =
\Big(\frac{2\alpha rp}{n^{1-r}}\Big)^{\frac{1}{p'}}<\Big(\frac{1}{7}\Big)^{\frac{1}{p'}}<1,
\end{equation}
and for  $ p'=1$ and $n\geq n_2(\alpha,r,p)$, taking into account decreasing on the interval  $[e, \infty)$ of the function $\frac{\ln \upsilon}{\upsilon}$, we have
$$
 \frac{2\alpha r}{n^{1-r}}J_{1}\Big(\frac{ \pi n^{1-r}}{\alpha r}\Big)=
 \frac{2\alpha r}{n^{1-r}}\int\limits_{0}^{\frac{\pi n^{1-r}}{\alpha r}}\frac{dt}{\sqrt{t^2+1}}<
 \frac{2\alpha r}{n^{1-r}}\bigg(1+\int\limits_{1}^{\frac{\pi n^{1-r}}{\alpha r}}\frac{dt}{\sqrt{t^2+1}}\bigg)<
 $$
\begin{equation}\label{in2}
<
\frac{2\alpha r}{n^{1-r}}+  \frac{2\alpha r}{n^{1-r}}\ln \Big(\frac{ \pi n^{1-r}}{\alpha r}\Big)\leq
    \frac{2}{14}+\frac{2\pi\ln14\pi}{14\pi}<1.
\end{equation}

Formulas (\ref{in1}) and (\ref{in2}) prove (\ref{f201}).
For $1\leq p'<\infty$ estimate (\ref{f109}) follows from  (\ref{f200}) and (\ref{f201}).

Let us verify validity of the estimate  (\ref{f109}) for  $p'=\infty$. It follows from  (\ref{r11}) and (\ref{xi}) that
\begin{equation}\label{z1}
 \big\|\mathcal{P}_{\alpha,r,n} \big\|_{\infty}= \sum\limits_{k=0}^{\infty}e^{-\alpha (k+n)^{r}}=
 \int\limits_{n}^{\infty}e^{-\alpha t^{r}}dt+\Theta_{\alpha,r,n}^{(11)}e^{-\alpha n^{r}}, \ |\Theta_{\alpha,r,n}^{(11)}|\leq1.
  \end{equation}

 Setting in formula  (\ref{0lemma_4}) $\gamma=\alpha$, $\delta=0$ and $m=n$, from (\ref{z1}) we obtain that for arbitrary  $n\geq n_2(\alpha,r,p)$
 \begin{equation}\label{z2}
 \big\|\mathcal{P}_{\alpha,r,n} \big\|_{\infty}=\frac{e^{-\alpha n^{r}}}{\alpha r}n^{1-r}\Big(
1+\Theta_{\alpha,r,n}^{(12)}\Big(\frac{1-r}{\alpha r}\frac{1}{n^{r}}+\frac{\alpha r}{n^{1-r}}\Big)\Big), \ \
 \end{equation}
where $\ |\Theta_{\alpha,r,n}^{(12)}|\leq\frac{14}{13}$.

For $p'=\infty$ the validity of (\ref{f109}) follows from  (\ref{z2}) and the equality  $J_{\infty}\big(\frac{\pi n^{1-r}}{\alpha r}\big)=1$.

To complete the proof of theorem 1 it suffices to find the upper estimate of the quantity   $M_{n}$ in formula (\ref{f4440}).
 It is clear that
 $$
M_{n}=
 \sup\limits_{t\in\mathbb{R}}
 \frac{\big|\mathcal{P'}_{\alpha,r,n}(t) \big|\big|\mathcal{P}_{\alpha,r,n}(t) \big|}
 {\big|\mathcal{P}_{\alpha,r,n}(t) \big|^{2}}=
 $$
\begin{equation}\label{ss4}
= \max\bigg\{\sup\limits_{|t|\leq\frac{\alpha r}{n^{1-r}}}
 \frac{\big|\mathcal{P'}_{\alpha,r,n}(t) \big|\big|\mathcal{P}_{\alpha,r,n}(t) \big|}
 {\big|\mathcal{P}_{\alpha,r,n}(t) \big|^{2}},
 \sup\limits_{\frac{\alpha r}{n^{1-r}}\leq|t|\leq\pi}
 \frac{\big|\mathcal{P'}_{\alpha,r,n}(t) \big|\big|\mathcal{P}_{\alpha,r,n}(t) \big|}
 {\big|\mathcal{P}_{\alpha,r,n}(t) \big|^{2}}\bigg\}.
\end{equation}

In view of formulas  (\ref{n1}) and (\ref{f34}) and the fact that $R_n(t)>0$
  for $n\geq n_2(\alpha,r,p)$ we obtain
\begin{equation}\label{0ss}
  \big|\mathcal{P}_{\alpha,r,n}(t) \big|^{2}>Q_{n}(t)>\frac{9}{14}\frac{ e^{-2\alpha n^{r}}}{t^{2}+(\alpha r n^{r-1})^{2}}.
\end{equation}

%Із (\ref{mmm1}) і (\ref{0ss}) отримуємо, що при $n\geq n_2$

It directly follows from  (\ref{r11}) that
\begin{equation}\label{p1}
   \big|\mathcal{P}_{\alpha,r,n}(t) \big| \leq
\sum\limits_{k=0}^{\infty}e^{-\alpha(k+n)^{r}},  \ \
\big|\mathcal{P'}_{\alpha,r,n}(t) \big| \leq
\sum\limits_{k=1}^{\infty}k e^{-\alpha(k+n)^{r}}.
\end{equation}

By virtue of (\ref{z2}) for $n\geq n_2(\alpha,r,p)$ we have
\begin{equation}\label{p3}
\big|\mathcal{P}_{\alpha,r,n}(t) \big| \leq\sum\limits_{k=0}^{\infty}e^{-\alpha (k+n)^{r}}<\frac{14}{13}e^{-\alpha n^{r}}\frac{n^{1-r}}{\alpha r}.
\end{equation}

The function $te^{-\alpha t^{r}}$
is monotone decreasing for $t>(\alpha r)^{-\frac{1}{r}}$. Therefore, according to  (\ref{xi}), for    ${n\geq n_2(\alpha,r,p)}$ the following estimate takes place
$$
\sum\limits_{k=1}^{\infty}e^{-\alpha (k+n)^{r}}k=\sum\limits_{k=n}^{\infty}e^{-\alpha k^{r}}k-
n\sum\limits_{k=n}^{\infty}e^{-\alpha k^{r}}\leq
$$
\begin{equation}\label{s1}
\leq e^{-\alpha n^{r}}n+\int\limits_{n}^{\infty}e^{-\alpha t^{r}}tdt-n\int\limits_{n}^{\infty}e^{-\alpha t^{r}}dt.
\end{equation}

Setting in   (\ref{0lemma_4}) $\gamma=\alpha$, $\delta=1$, $m=n$, and also    $\gamma=\alpha$, $\delta=0$, $m=n$, from  (\ref{p1}) and (\ref{s1}) we have
\begin{equation}\label{s2}
\big|\mathcal{P'}_{\alpha,r,n}(t) \big| \leq
e^{-\alpha n^{r}}\Big(\frac{42}{13}\Big(\frac{n^{1-r}}{\alpha r}\Big)^{2}+n\Big), \ \ n\geq n_2(\alpha,r,p).
\end{equation}

In view of (\ref{0ss}), (\ref{p3}) and (\ref{s2}) for   $n\geq n_2(\alpha,r,p)$ we arrive at the estimate
$$
\sup\limits_{|t|\leq\frac{\alpha r}{n^{1-r}}}
 \frac{\big|\mathcal{P'}_{\alpha,r,n}(t) \big|\big|\mathcal{P}_{\alpha,r,n}(t) \big|}
 {\big|\mathcal{P}_{\alpha,r,n}(t) \big|^{2}}\leq
$$
$$
\leq \frac{14}{9}e^{2\alpha n^{r}}
 \sup\limits_{|t|<\frac{\alpha r}{n^{1-r}}}
 \big|\mathcal{P'}_{\alpha,r,n}(t) \big|\big|\mathcal{P}_{\alpha,r,n}(t) \big|
 (t^{2}+\Big(\frac{\alpha r}{n^{1-r}}\Big)^{2})\leq
 $$
\begin{equation}\label{p6}
  \leq \frac{5488}{507}
\Big(\Big(\frac{n^{1-r}}{\alpha r}\Big)^{2}+n\Big)\frac{n^{1-r}}{\alpha r}
\Big(\frac{\alpha r}{n^{1-r}}\Big)^{2}=\frac{5488}{507}
\Big(\frac{n^{1-r}}{\alpha r}+ \alpha rn^{r}\Big).
\end{equation}

Applying  the Abel transformation to the function $\mathcal{P}_{\alpha,r,n}(t)$ for $0<|t|\leq\pi$,
and taking into account the inequallity
$$
\Big|\sum\limits_{j=0}^{k}e^{ijt}\Big|\leq\frac{\pi}{|t|}, \ 0<|t|\leq\pi,
$$
we get
\begin{equation}\label{g1}
\big|\mathcal{P}_{\alpha,r,n}(t) \big| =
\Big|
\sum\limits_{k=0}^{\infty}(e^{-\alpha (k+n)^{r}}-e^{-\alpha (k+n+1)^{r}})\sum\limits_{j=0}^{k}e^{ijt}\Big|
\leq\frac{\pi}{|t|}e^{-\alpha n^{r}}.
\end{equation}
By analogy, for $0<|t|\leq\pi$
$$
\big|\mathcal{P'}_{\alpha,r,n}(t) \big|
=\Big|
\sum\limits_{k=0}^{\infty}(e^{-\alpha (k+n)^{r}}k-e^{-\alpha (k+n+1)^{r}}(k+1))\sum\limits_{j=0}^{k}e^{ijt}\Big|\leq
$$
$$
\leq\frac{\pi}{|t|}\sum\limits_{k=0}^{\infty}|e^{-\alpha (k+n)^{r}}k-e^{-\alpha (k+n+1)^{r}}(k+1)|\leq
$$
\begin{equation}\label{11p7}
\leq \frac{\pi}{|t|}\bigg(\sum\limits_{k=0}^{\infty}k(e^{-\alpha (k+n)^{r}}-e^{-\alpha (k+n+1)^{r}})+
\sum\limits_{k=0}^{\infty}e^{-\alpha (k+n+1)^{r}}\bigg)=
\end{equation}
According to (\ref{p3}) and (\ref{11p7})
\begin{equation}\label{p7}
 \big|\mathcal{P'}_{\alpha,r,n}(t) \big|\leq
\frac{2\pi}{|t|}\sum\limits_{k=0}^{\infty}e^{-\alpha (k+n+1)^{r}}\leq \frac{28\pi}{13|t|}
e^{-\alpha n^{r}}\frac{n^{1-r}}{\alpha r}.
\end{equation}

In view of (\ref{0ss}), (\ref{g1}) and (\ref{p7}) we obtain the estimate
 $$
 \sup\limits_{\frac{\alpha r}{n^{1-r}}\leq|t|\leq\pi}
 \frac{\big|\mathcal{P'}_{\alpha,r,n}(t) \big|\big|\mathcal{P}_{\alpha,r,n}(t) \big|}
 {\big|\mathcal{P}_{\alpha,r,n}(t) \big|^{2}}\leq
 $$
 $$
\leq\frac{14}{9}
 e^{2\alpha n^{r}}
 \sup\limits_{\frac{\alpha r}{n^{1-r}}\leq|t|\leq\pi}
 \big|\mathcal{P'}_{\alpha,r,n}(t) \big|\big|\mathcal{P}_{\alpha,r,n}(t) \big|
 \Big(t^{2}+\Big(\frac{\alpha r}{n^{1-r}}\Big)^{2}\Big)\leq
$$
\begin{equation}\label{p8}
\leq\frac{392\pi^{2}}{117}\frac{n^{1-r}}{\alpha r}\sup\limits_{ \frac{\alpha r}{n^{1-r}}\leq|t|\leq\pi}
\frac{t^{2}+\big(\frac{\alpha r}{n^{1-r}}\big)^{2}}{t^{2}}\leq\frac{784\pi^{2}}{117}\frac{n^{1-r}}{\alpha r}.
\end{equation}

Combining  (\ref{ss4}), (\ref{p6}) and (\ref{p8}), we arrive at the eatimate
\begin{equation}\label{p9}
 M_{n}\leq\frac{784\pi^{2}}{117}\Big(\frac{n^{1-r}}{\alpha r}+ \alpha rn^{r}\Big), \ \ n\geq n_2(\alpha,r,p).
\end{equation}

It follows from conditions  (\ref{n_p}) and (\ref{n1}) that  $n_0(\alpha,r,p)\geq n_2(\alpha,r,p)$ for arbitrary ${1\leq p\leq\infty}$. It means that estimates (\ref{f109}) and (\ref{p9}) are true also for   ${n\geq n_0(\alpha,r,p)}$.
Let us show that for  $n\geq n_0(\alpha,r,p)$ the condition (\ref{n_nom}) is satisfied. This is obvious for $p'=\infty$. For
 $1\leq p'<\infty$ by virtue of  (\ref{p9}), we have
 \begin{equation}\label{c11}
  4\pi M_{n}p'\leq\frac{3136\pi^{3}}{117}\Big(\frac{n^{1-r}}{\alpha r}+ \alpha rn^{r}\Big)p'<
  27\pi^{3}\Big(\frac{n^{1-r}}{\alpha r}+ \alpha r\chi(p)n^{r}\Big)p'.
 \end{equation}
 According to (\ref{n_p}) and (\ref{c11})
 for any
 $ n\geq n_{0}(\alpha,r,p)$
the following inequality is true
$$
4\pi p'M_{n}\leq n,
$$
which is equivalent to  (\ref{n_nom}) for $1\leq p'<\infty$.

By using formulas (\ref{f4440}),  (\ref{f109}) and (\ref{p9}) for $n\geq n_0(\alpha,r,p)$   we arrive at the estimate
$$
{\cal E}_{n}(C^{\alpha,r}_{\beta,p})_{C}=
$$
$$
=e^{-\alpha n^{r}}n^{\frac{1-r}{p}}\bigg(\frac{2^{\frac{1}{p'}}}{(\alpha r)^{\frac{1}{p}}}J_{p'}\Big(\frac{ \pi n^{1-r}}{\alpha r}\Big)+\Theta_{\alpha,r,p,n}^{(3)}\Big(\frac{1-r}{(\alpha r)^{1+\frac{1}{p}}}J_{p'}\Big(\frac{ \pi n^{1-r}}{\alpha r}\Big)\frac{1}{n^{r}}+\frac{1}{n^{\frac{1-r}{p}}}\Big)\bigg)\times
$$
\begin{equation}\label{subst}
 \times\Big(\frac{\|\cos t\|_{p'}}{2^{\frac{1}{p'}}\pi^{1+\frac{1}{p'}}}
+\delta_{n}^{(3)}\Big(\frac{1}{\alpha r}\frac{1}{n^{r}}+ \frac{\alpha r}{n^{1-r}}\Big)\Big), \ 1\leq p\leq\infty,
 \end{equation}
where for  $\Theta_{\alpha,r,p,n}^{(3)}$ the estimate (\ref{theta3}) takes place, and $|\delta_{n}^{(3)}|<\frac{10976\pi^{2}}{117}$.

For $n\geq n_0(\alpha,r,p)$ the following inequality holds
$$
|\delta_{n}^{(3)}|\frac{2^{\frac{1}{p'}}}{(\alpha r)^{\frac{1}{p}}}J_{p'}\Big(\frac{ \pi n^{1-r}}{\alpha r}\Big)\Big(\frac{1}{\alpha r}\frac{1}{n^{r}}+ \frac{\alpha r}{n^{1-r}}\Big)<
$$
\begin{equation}\label{subst1}
<\frac{21952\pi^2}{117}
\Big(\frac{1}{(\alpha r)^{1+\frac{1}{p}}}J_{p'}\Big(\frac{ \pi n^{1-r}}{\alpha r}\Big)\frac{1}{n^{r}}+\frac{1}{n^{\frac{1-r}{p}}}\Big),
\end{equation}
which follows from  (\ref{f201}) for $1\leq p'<\infty$, and it is obvious for $p'=\infty$. Besides, according to  (\ref{theta3}) and (\ref{n_p}) for $n\geq n_0(\alpha,r,p)$
$$
\Big|\Theta_{\alpha,r,p,n}^{(3)}\Big|\Big(\frac{1-r}{(\alpha r)^{1+\frac{1}{p}}}J_{p'}\Big(\frac{ \pi n^{1-r}}{\alpha r}\Big)\frac{1}{n^{r}}+\frac{1}{n^{\frac{1-r}{p}}}\Big)
 \Big(\frac{\|\cos t\|_{p'}}{2^{\frac{1}{p'}}\pi^{1+\frac{1}{p'}}}
+|\delta_{n}^{(3)}|\Big(\frac{1}{\alpha r}\frac{1}{n^{r}}+ \frac{\alpha r}{n^{1-r}}\Big)\Big)<
$$
\begin{equation}\label{subst2}
<\frac{363\pi^{2}}{50}\Big(\frac{1-r}{(\alpha r)^{1+\frac{1}{p}}}J_{p'}\Big(\frac{ \pi n^{1-r}}{\alpha r}\Big)\frac{1}{n^{r}}+\frac{1}{n^{\frac{1-r}{p}}}\Big).
\end{equation}

In view of formulas  (\ref{subst})--(\ref{subst2}) we arrive at (\ref{theorem1}).
Theorem 1 is proved.

\vskip 10mm

{ \bf 4. Proof of lemma $1$.} It is obvious that for  $ 1\leq s\leq \infty$
$$
\inf\limits_{\lambda\in\mathbb{R}}\|\phi(t)-\lambda\|_{s}\leq \|\phi\|_{s}, \
$$
$$
\frac{1}{2}\|\phi\big(t+\frac{\pi}{n}\big)-\phi(t)\|_{s}\leq
\sup\limits_{h\in\mathbb{R}}\frac{1}{2}\|\phi\big(t+h\big)-\phi(t)\|_{s}
$$
and
$$
\sup\limits_{h\in\mathbb{R}}\frac{1}{2}\|\phi\big(t+h\big)-\phi(t)\|_{s}\leq
\inf\limits_{\lambda\in\mathbb{R}}\|\phi(t)-\lambda\|_{s}.
$$
Hence, in order to proof lemma it suffices to verify the validity of formula (\ref{a1}) and relation
\begin{equation}\label{a2}
 \frac{1}{2}\|\phi\big(t+\frac{\pi}{n}\big)-\phi(t)\|_{s}
 \geq\|r\|_{s}\Big(\frac{\|\cos t\|_{s}}{(2\pi)^{\frac{1}{s}}}-14\pi\frac{M}{n}\Big).
\end{equation}

First, we consider the case  $1\leq s<\infty$.  Let verify the validity of equality  (\ref{a1}). Setting
\begin{equation}\label{phi_k}
 \phi_{k}(t)=g\Big(\frac{k\pi}{n}\Big)\cos(nt+\gamma)
 +h\Big(\frac{k\pi}{n}\Big)\sin(nt+\gamma), \ \ k=\overline{-n+1, \ n},
\end{equation}
we get
$$
 \|\phi\|_{s}=\bigg(\sum\limits_{k=-n+1}^{n}\int\limits_{\frac{(k-1)\pi}{n}}^{\frac{k\pi}{n}}|\phi(t)|^{s}dt\bigg)^{\frac{1}{s}}=
$$
\begin{equation}\label{a3}
=\bigg(\sum\limits_{k=-n+1}^{n}\int\limits_{\frac{(k-1)\pi}{n}}^{\frac{k\pi}{n}}|\phi_{k}(t)|^{s}dt\bigg)^{\frac{1}{s}}+
\Theta_{n}^{(1)}\bigg(\sum\limits_{k=-n+1}^{n}\int\limits_{\frac{(k-1)\pi}{n}}^{\frac{k\pi}{n}}|\phi(t)-\phi_{k}(t)|^{s}dt\bigg)^{\frac{1}{s}}, \ \ |\Theta_{n}^{(1)}|\leq1.
\end{equation}
Let us find the estimate of first term in  (\ref{a3}). It is obvious, that according to  (\ref{phi_k})
$$
\bigg(\sum\limits_{k=-n+1}^{n}\int\limits_{\frac{(k-1)\pi}{n}}^{\frac{k\pi}{n}}|\phi_{k}(t)|^{s}dt\bigg)^{\frac{1}{s}}=
$$
$$
=\bigg(\sum\limits_{k=-n+1}^{n}r^{s}\Big(\frac{k\pi}{n}\Big)\int\limits_{\frac{(k-1)\pi}{n}}^{\frac{k\pi}{n}}
\Big|\cos \Big(nt+\gamma-\mathrm{arg}\Big(g\Big(\frac{k\pi}{n}\Big)+ih\Big(\frac{k\pi}{n}\Big) \Big)\Big)\Big|^{s}dt\bigg)^{\frac{1}{s}}=
$$
\begin{equation}\label{a21}
=\bigg(\sum\limits_{k=-n+1}^{n}r^{s}\Big(\frac{k\pi}{n}\Big)\frac{1}{n}\int\limits_{0}^{\pi}|\cos t|^{s}dt\bigg)^{\frac{1}{s}}=
\frac{\|\cos t\|_{s}}{(2\pi)^{\frac{1}{s}}}\bigg(\sum\limits_{k=-n+1}^{n}r^{s}\Big(\frac{k\pi}{n}\Big)\frac{\pi}{n}\bigg)^{\frac{1}{s}},
\end{equation}
where
$r(t)$ is defined by formula (\ref{r}), and $i$ is imaginary unit.

Let us show that for any collection of points  $\xi_{k}, \ k=\overline{-n+1, \ n}$, such that ${\frac{(k-1)\pi}{n}\leq \xi_{k}\leq\frac{k\pi}{n}}$,
for $n\geq4\pi s M$
 the following estimate is true
\begin{equation}\label{sys}
  \bigg(\sum\limits_{k=-n+1}^{n}r^{s}(\xi_{k})\frac{\pi}{n}\bigg)^{\frac{1}{s}}=
  \|r\|_{s}\Big(1+\Theta_{n}^{(2)}\frac{M}{n} \Big), \ \ |\Theta_{n}^{(2)}|\leq4.
\end{equation}
%де $\delta_{2}=\delta_{2}(\xi;s;n):=\frac{n}{M \|r\|_{s}}\bigg(\sum\limits_{k=-n+1}^{n}r^{s}(\xi_{k})\frac{\pi}{n}\bigg)^{\frac{1}{s}}-\frac{n}{M}$.

Indeed, since
$$
\sum\limits_{k=-n+1}^{n}r^{s}(\xi_{k})\frac{\pi}{n}=\int\limits_{-\pi}^{\pi}r^{s}(t)dt+
\Theta_{n}^{(3)}\frac{\bigvee\limits_{-\pi}^{\pi}(r^{s})}{n}, \ |\Theta_{n}^{(3)}|\leq\pi,
$$
 and under the condition
 \begin{equation}\label{condit}
n\geq \frac{2\pi\bigvee\limits_{-\pi}^{\pi}(r^{s})}{\|r\|_{s}^{s}}
\end{equation}
\begin{equation}\label{a22}
 \bigg( \int\limits_{-\pi}^{\pi}r^{s}(t)dt+
\Theta_{n}^{(3)}\frac{\bigvee\limits_{-\pi}^{\pi}(r^{s})}{n}\bigg)^{\frac{1}{s}}=
\|r\|_{s}\bigg(1+\Theta_{n}^{(4)}\frac{\bigvee\limits_{-\pi}^{\pi}(r^{s})}{ns\|r\|_{s}^{s}}\bigg), \ \ |\Theta_{n}^{(4)}|\leq2,
\end{equation}
hence
\begin{equation}\label{a4}
\bigg(\sum\limits_{k=-n+1}^{n}r^{s}(\xi_{k})\frac{\pi}{n}\bigg)^{\frac{1}{s}}
 =\|r\|_{s}\bigg(1+\Theta_{n}^{(4)}\frac{\bigvee\limits_{-\pi}^{\pi}(r^{s})}{ns\|r\|_{s}^{s}}\bigg), \ \ |\Theta_{n}^{(4)}|\leq2.
\end{equation}

%Із (\ref{a4})  та співвідношень
It is easy to verify that
\begin{equation}\label{var}
 \bigvee\limits^{\pi}_{-\pi}(r^{s})=s\int\limits_{-\pi}^{\pi}r^{s-1}(t)|r'(t)|dt\leq s\|r\|_{s}^{s} \Big\|\frac{r'(t)}{r(t)}\Big\|_{\infty},
\end{equation}
\begin{equation}\label{a70}
 \Big|\frac{r'(t)}{r(t)}\Big|= \Big|\frac{g(t)g'(t)+h(t)h'(t)}{r^{2}(t)}\Big|
 \leq\frac{|g'(t)|+|h'(t)|}{r(t)}\leq 2M, \ \ t\in\mathbb{R},
\end{equation}
therefore
\begin{equation}\label{variation}
 \frac{\bigvee\limits_{-\pi}^{\pi}(r^{s})}{\|r\|_{s}^{s}}\leq   s\Big\|\frac{r'(t)}{r(t)}\Big\|_{\infty}\leq 2 s M.
\end{equation}
By virtue of   (\ref{variation}), for $n\geq 4\pi s M$ the condition (\ref{condit}) is satisfied. Therefore, according to (\ref{a4}),
the estimate (\ref{sys}) takes place. Setting in   (\ref{sys})  $\xi_{k}=\frac{k\pi}{n}$, $k=\overline{-n+1,n}$,  in view of (\ref{a21})
we obtain
\begin{equation}\label{a6}
  \bigg(\sum\limits_{k=-n+1}^{n}\int\limits_{\frac{(k-1)\pi}{n}}^{\frac{k\pi}{n}}|\phi_{k}(t)|^{s}dt\bigg)^{\frac{1}{s}}=
\|r\|_{s}\Big(\frac{\|\cos t\|_{s}}{(2\pi)^{\frac{1}{s}}}+\Theta_{n}^{(5)}\frac{M}{n} \Big), \ \ |\Theta_{n}^{(5)}|\leq4.
\end{equation}
Let us find upper estimate of the second term in   (\ref{a3}). On the basis of   (\ref{snow}) and (\ref{phi_k})
$$
\phi(t)-\phi_{k}(t)=
$$
$$
=\Big(r(t)-r\Big(\frac{k\pi}{n}\Big)\Big)\Big(\frac{g(\frac{k\pi}{n})}{r(\frac{k\pi}{n})}\cos(nt+\gamma)+
\frac{h(\frac{k\pi}{n})}{r(\frac{k\pi}{n})}\sin(nt+\gamma)\Big)+
$$
\begin{equation}\label{a7}
+r(t)\Big(\Big(\frac{g(t)}{r(t)}-
\frac{g(\frac{k\pi}{n})}{r(\frac{k\pi}{n})}\Big)\cos(nt+\gamma)+ \Big
(\frac{h(t)}{r(t)}-
\frac{h(\frac{k\pi}{n})}{r(\frac{k\pi}{n})}\Big)\sin(nt+\gamma)\Big),
\end{equation}
therefore
\begin{equation}\label{a9}
\bigg(\sum\limits_{k=-n+1}^{n}\int\limits_{\frac{(k-1)\pi}{n}}^{\frac{k\pi}{n}}|\phi(t)-\phi_{k}(t)|^{s}dt\bigg)^{\frac{1}{s}}\leq I_{n}^{(1)}+ I_{n}^{(2)},
\end{equation}
where
$$
 I_{n}^{(1)}:=\bigg(\sum\limits_{k=-n+1}^{n}\int\limits_{\frac{(k-1)\pi}{n}}^{\frac{k\pi}{n}}
\Big|r(t)-r\Big(\frac{k\pi}{n}\Big)\Big|^{s}(|\cos(nt+\gamma)|+|\sin(nt+\gamma)|)^{s}dt\bigg)^{\frac{1}{s}},
$$
$$
I_{n}^{(2)}:=\bigg(\sum\limits_{k=-n+1}^{n}\int\limits_{\frac{(k-1)\pi}{n}}^{\frac{k\pi}{n}}
r^{s}(t)\Big(\Big|\frac{g(t)}{r(t)}-
\frac{g(\frac{k\pi}{n})}{r(\frac{k\pi}{n})}\Big||\cos(nt+\gamma)|+
$$
$$
+\Big
|\frac{h(t)}{r(t)}-
\frac{h(\frac{k\pi}{n})}{r(\frac{k\pi}{n})}\Big||\sin(nt+\gamma)|\Big)^{s}dt\bigg)^{\frac{1}{s}}.
$$

 Using obvious inequality
 \begin{equation}\label{cos_sin}
  |\cos t|+|\sin t|\leq \sqrt{2},
\end{equation}
Lagrange theorem and relation  (\ref{a70}),  we have
$$
 I_{n}^{(1)}\leq\sqrt{2}\bigg(\sum\limits_{k=-n+1}^{n}
 \max\limits_{\frac{(k-1)\pi}{n}\leq t \leq \frac{k\pi}{n}} \Big|r(t)-r\Big(\frac{k\pi}{n}\Big)\Big|^{s}\frac{\pi}{n}
\bigg)^{\frac{1}{s}}\leq
$$
\begin{equation}\label{a30}
  \leq\frac{\sqrt{2}\pi}{n}\sup\limits_{t\in\mathbb{R}}\Big|\frac{r'(t)}{r(t)}\Big|
\bigg(\sum\limits_{k=-n+1}^{n}
  \max\limits_{\frac{(k-1)\pi}{n}\leq t \leq \frac{k\pi}{n}}r^{s}(t)\frac{\pi}{n}
\bigg)^{\frac{1}{s}}.
\end{equation}

It is follows from  (\ref{sys}), (\ref{a70}) and (\ref{a30}),  that for $n\geq4\pi s M$
\begin{equation}\label{a32}
  I_{n}^{(1)}\leq2\sqrt{2}\pi\frac{M}{n}(1+4\frac{M}{n})\|r\|_{s}\leq
2\sqrt{2}\pi\frac{M}{n}\Big(1+\frac{1}{\pi}\Big)\|r\|_{s}=\frac{2\sqrt{2}M(1+\pi)}{n}\|r\|_{s}.
\end{equation}

It is easy to see that
\begin{equation}\label{aa13}
I_{n}^{(2)}\leq\bigg(\sum\limits_{k=-n+1}^{n}
\Big(\max\limits_{\frac{(k-1)\pi}{n}\leq t \leq \frac{k\pi}{n}}\Big\{\Big|\frac{g(t)}{r(t)}-
\frac{g(\frac{k\pi}{n})}{r(\frac{k\pi}{n})}\Big||\cos(nt+\gamma)|+$$
$$
+\Big
|\frac{h(t)}{r(t)}-
\frac{h(\frac{k\pi}{n})}{r(\frac{k\pi}{n})}\Big
||\sin(nt+\gamma)|\Big\}\Big)^{s}
\int\limits_{\frac{(k-1)\pi}{n}}^{\frac{k\pi}{n}}
r^{s}(t)dt\bigg)^{\frac{1}{s}}.
\end{equation}

For any  $t_{1}, \ t_{2}\in\mathbb{R}$ such that $|t_{1}-t_{2}|\leq\frac{\pi}{n}$ the following inequalities take place
\begin{equation}\label{t1}
  \Big|\frac{g(t_{1})}{r(t_{1})}-
\frac{g(t_{2})}{r(t_{2})}\Big|
 \leq
 \frac{3\pi M}{n},
\end{equation}
\begin{equation}\label{t2}
   \Big|\frac{h(t_{1})}{r(t_{1})}-
\frac{h(t_{2})}{r(t_{2})}\Big| \leq
 \frac{3\pi M}{n}.
\end{equation}

Indeed, by virtue of Lagrange theorem, taking into account  (\ref{M}) and (\ref{a70}),  we have
$$
 \Big|\frac{g(t_{1})}{r(t_{1})}-
\frac{g(t_{2})}{r(t_{2})}\Big|\leq
 \frac{\pi}{n}
\sup\limits_{t\in\mathbb{R}}\Big|\frac{g'(t)r(t)-g(t)r'(t)}{r^{2}(t)}\Big|
\leq
$$
\begin{equation}\label{a500}
\leq\frac{\pi}{n}
\sup\limits_{t\in\mathbb{R}}\frac{|g'(t)|}{r(t)}+\frac{\pi}{n}\sup\limits_{t\in\mathbb{R}}
\frac{|r'(t)|}{r(t)}
\leq\frac{3\pi M}{n}.
\end{equation}
By analogy, we prove the inequality   (\ref{t2}).
In view of (\ref{cos_sin}), (\ref{t1}), (\ref{t2}) and (\ref{aa13}) we obtain
\begin{equation}\label{a13}
I_{n}^{(2)}
\leq\frac{3\sqrt{2}\pi M}{n}
\|r\|_{s}, \ \ n\in\mathbb{N}.
\end{equation}

Combining  (\ref{a9}), (\ref{a32}) and (\ref{a13}),  we arrive at the estimate
\begin{equation}\label{a14}
\bigg(\sum\limits_{k=-n+1}^{n}\int\limits_{\frac{(k-1)\pi}{n}}^{\frac{k\pi}{n}}|\phi(t)-\phi_{k}(t)|^{s}dt\bigg)^{\frac{1}{s}}
\leq\sqrt{2}(5\pi+2)\|r\|_{s}\frac{M}{n}, \ \ n\geq4\pi s M.
\end{equation}

By comparing estimates  (\ref{a3}), (\ref{a6}) and (\ref{a14}) we conclude that for   ${n\geq4\pi s M}$
\begin{equation}\label{abc2}
\|\phi\|_{s}= \|r\|_{s}\Big(\frac{\|\cos t\|_{s}}{(2\pi)^{\frac{1}{s}}}+\delta_{s,n}^{(1)}\frac{M}{n}\Big), \ |\delta_{s,n}^{(1)}|\leq \sqrt{2}(5\pi+2)+4, \ 1\leq s<\infty.
\end{equation}

Further, we prove the relation  (\ref{a2}) for $1\leq s<\infty$. In view of definition  (\ref{snow})
$$
|\phi\big(t+\frac{\pi}{n}\big)-\phi(t)|=
$$
$$
=\Big|2\phi(t)+
g\Big(t+\frac{\pi}{n}\Big)\cos(nt+\gamma )
+h\Big(t+\frac{\pi}{n}\Big)\sin(nt+\gamma )-
$$
$$
-(g(t)\cos(nt+\gamma )
+h(t)\sin(nt+\gamma))\Big|=
$$
$$
 =\Big|2\phi(t)+\Big(r\Big(t+\frac{\pi}{n}\Big)-r(t)\Big)\Big(
 \frac{g(t+\frac{\pi}{n})}{r(t+\frac{\pi}{n})}\cos(nt+\gamma )
+\frac{h(t+\frac{\pi}{n})}{r(t+\frac{\pi}{n})}\sin(nt+\gamma)\Big)+
$$
\begin{equation}\label{a15}
+r(t)\Big(\Big(
 \frac{g(t+\frac{\pi}{n})}{r(t+\frac{\pi}{n})}-\frac{g(t)}{r(t)}\Big)\cos(nt+\gamma )
+\Big(\frac{h(t+\frac{\pi}{n})}{r(t+\frac{\pi}{n})}-\frac{h(t)}{r(t)}\Big)
\sin(nt+\gamma)\Big)\Big|,
\end{equation}
therefore for any  $1\leq s\leq\infty$ by virtue of  (\ref{cos_sin}), (\ref{t1}) and (\ref{t2}), we get
$$
\frac{1}{2}\|\phi\big(t+\frac{\pi}{n}\big)-\phi(t)\|_{s}\geq
$$
\begin{equation}\label{a16}
\geq
\|\phi\|_{s}-\frac{1}{\sqrt{2}}\bigg(\Big\| r\Big(t+\frac{\pi}{n}\Big)-r(t)\Big \|_{s}
+3\pi\|r\|_{s}\frac{M}{n}\bigg).
\end{equation}
By applying the Lagrange theorem, we obtain
$$
\Big\| r\Big(t+\frac{\pi}{n}\Big)-r(t)\Big \|_{s}
  =\bigg(\sum\limits_{k=-n+1}^{n}\int\limits_{\frac{(k-1)\pi}{n}}^{\frac{k\pi}{n}}
\Big|r\Big(t+\frac{\pi}{n}\Big)-r(t)\Big|^{s}dt\bigg)^{\frac{1}{s}}\leq
$$
$$
 \leq\bigg(\sum\limits_{k=-n+1}^{n}
 \max\limits_{\frac{(k-1)\pi}{n}\leq t \leq \frac{k\pi}{n}} \Big|r\Big(t+\frac{\pi}{n}\Big)-r(t)  \Big|^{s}\frac{\pi}{n}
\bigg)^{\frac{1}{s}}\leq
$$
\begin{equation}\label{a33}
  \leq\frac{\pi}{n}\sup\limits_{t\in\mathbb{R}}\Big|\frac{r'(t)}{r(t)}\Big|
\bigg(\sum\limits_{k=-n+1}^{n}
 \max\limits_{\frac{(k-1)\pi}{n}\leq t \leq \frac{k\pi}{n}}r^{s}(t)\frac{\pi}{n}
\bigg)^{\frac{1}{s}}, \ 1\leq s<\infty.
\end{equation}

It follows from   (\ref{sys}), (\ref{a70}) and (\ref{a33}) that for $n\geq4\pi s M$
\begin{equation}\label{a18}
\Big\|r\Big(t+\frac{\pi}{n}\Big)-r(t)\Big \|_{s}\leq
(2\pi+2)\|r \|_{s}\frac{M}{n}.
\end{equation}

In view of   (\ref{abc2}),  (\ref{a16}) and (\ref{a18})  for  $n\geq 4\pi s M$ we arrive at the estimate
  $$
\frac{1}{2}\|\phi\big(t+\frac{\pi}{n}\big)-\phi(t)\|_{s}\geq\|\phi\|_{s}
-\frac{5\pi+2}{\sqrt{2}}\|r \|_{s}\frac{M}{n}\geq
$$
$$
\geq\|r\|_{s}\Big(\frac{\|\cos t\|_{s}}{(2\pi)^{\frac{1}{s}}}-\Big(\frac{15\pi+6}{\sqrt{2}}+4\Big)\frac{M}{n}\Big)>\|r\|_{s}\Big(\frac{\|\cos t\|_{s}}{(2\pi)^{\frac{1}{s}}}-14\pi\frac{M}{n}\Big), \ 1\leq s<\infty.
$$
 Thus, the validity of formula  (\ref{a2}) is established for $1\leq s<\infty$.

 Let us prove the relation (\ref{a1}) for $s=\infty$. Consider a function  $\phi^{*}(t)$ such that
 $$
 \phi^{*}(t)=\phi^{*}_{k}(t), \  \ \frac{(k-1)\pi}{n}\leq t \leq \frac{k\pi}{n}, \ k=\overline{-n+1, \ n},
 $$
 where
 \begin{equation}\label{b20}
  \phi^{*}_{k}(t)=g(t_{k}^{*})\cos (nt+\gamma)+h(t_{k}^{*})\sin (nt+\gamma),
 \end{equation}
  and points $t_{k}^{*}, \ t_{k}^{*}\in[\frac{(k-1)\pi}{n}, \frac{k\pi}{n}]$ are chosen from the condition
 $$
 r(t_{k}^{*})=\max\limits_{\frac{(k-1)\pi}{n}\leq t \leq \frac{k\pi}{n}}r(t).
 $$
 For the function  $\phi^{*}(t)$  the following equality takes place
 \begin{equation}\label{b21}
 \|\phi^{*}\|_{\infty}=\|r\|_{C}.
 \end{equation}

 Indeed,
  $$
 \|\phi^{*}\|_{\infty}=\max\limits_{-n+1\leq k\leq n}\mathop{\rm{ess}\sup}\limits_{\frac{(k-1)\pi}{n}\leq t \leq \frac{k\pi}{n}}|\phi^{*}(t)|=
  $$
  $$
  =\max\limits_{-n+1\leq k\leq n}r\big(t_{k}^{*}\big)\mathop{\max}\limits_{\frac{(k-1)\pi}{n}\leq t \leq \frac{k\pi}{n}}\Big|\frac{g(t_{k}^{*})}{r(t_{k}^{*})}\cos(nt+\gamma)+
\frac{h(t_{k}^{*})}{r(t_{k}^{*})}\sin(nt+\gamma)\Big|=
  $$
  $$
  =\max\limits_{-n+1\leq k\leq n}r\big(t_{k}^{*}\big)\mathop{\max}\limits_{\frac{(k-1)\pi}{n}\leq t \leq \frac{k\pi}{n}}
  \Big|\cos\Big(nt+\gamma-arg(g(t_{k}^{*})+ih(t_{k}^{*}))\Big)\Big|=
  $$
  $$
  =\max\limits_{-n+1\leq k\leq n}r\big(t_{k}^{*}\big) \ \|\cos t\|_{C}=
    \|r\|_{C}.
  $$

 It is obvious that in view of  (\ref{b21}) we obtain
 \begin{equation}\label{b1}
   \|\phi\|_{\infty}=\|\phi^{*}\|_{\infty}+\Theta_{n}^{(6)} \|\phi-\phi^{*}\|_{\infty}=
   \|r\|_{C}+\Theta_{n}^{(6)} \|\phi-\phi^{*}\|_{\infty}, \ \ |\Theta_{n}^{(6)}|\leq1.
 \end{equation}

  Let us find upper estimate for the quantity  $\|\phi-\phi^{*}\|_{\infty}$.
By virtue of  (\ref{snow}) and (\ref{b20}), for any $t\in[\frac{(k-1)\pi}{n}, \frac{k\pi}{n}]$ the following equality takes place
  $$
|\phi(t)-\phi_{k}^{*}(t)|=\Big|(r(t)-r(t_{k}^{*}))\Big(\frac{g(t_{k}^{*})}{r(t_{k}^{*})}\cos(nt+\gamma)+
\frac{h(t_{k}^{*})}{r(t_{k}^{*}))}\sin(nt+\gamma)\Big)+
$$
  \begin{equation}\label{b7}
+r(t)\Big(\Big(\frac{g(t)}{r(t)}-
\frac{g(t_{k}^{*})}{r(t_{k}^{*})}\Big)\cos(nt+\gamma)+ \Big
(\frac{h(t)}{r(t)}-
\frac{h(t_{k}^{*})}{r(t_{k}^{*})}\Big)\sin(nt+\gamma)\Big)\Big|.
\end{equation}

 By using (\ref{cos_sin}), the Lagrange theorem and inequality (\ref{a70}), we get
 $$
  \mathop{\rm{ess}\sup}\limits_{\frac{(k-1)\pi}{n}\leq t \leq \frac{k\pi}{n}}\Big|(r(t)-r(t_{k}^{*}))\Big(\frac{g(t_{k}^{*})}{r(t_{k}^{*})}\cos(nt+\gamma)+
\frac{h(t_{k}^{*})}{r(t_{k}^{*}))}\sin(nt+\gamma)\Big|\leq
 $$
 \begin{equation}\label{b8}
 \leq\sqrt{2} \mathop{\rm{ess}\sup}\limits_{\frac{(k-1)\pi}{n}\leq t \leq \frac{k\pi}{n}}\Big|r(t)-r(t_{k}^{*})\Big|\leq
  \frac{\sqrt{2}\pi}{n}\sup\limits_{t\in\mathbb{R}}\Big|\frac{r'(t)}{r(t)}\Big|\|r\|_{C}
  \leq
  \frac{2\sqrt{2}\pi M}{n}\|r\|_{C}.
 \end{equation}

 Further, it follows from (\ref{cos_sin}), (\ref{t1}) and (\ref{t2}) that
 $$
  \mathop{\rm{ess}\sup}\limits_{\frac{(k-1)\pi}{n}\leq t \leq \frac{k\pi}{n}}
 r(t)\Big(\Big|\frac{g(t)}{r(t)}-
\frac{g(t_{k}^{*})}{r(t_{k}^{*})}\Big||\cos(nt+\gamma)|+ \Big|\frac{h(t)}{r(t)}-
\frac{h(t_{k}^{*})}{r(t_{k}^{*})}\Big||\sin(nt+\gamma)|\Big) \leq
 $$
  \begin{equation}\label{b9}
 \leq3\sqrt{2}\pi
\frac{M}{n}\|r\|_{C}.
 \end{equation}
In view of (\ref{b7})--(\ref{b9}), we arrive at the estimate
 \begin{equation}\label{b10}
   \|\phi-\phi^{*}\|_{\infty}=\max\limits_{-n+1\leq k\leq n}\mathop{\rm{ess}\sup}\limits_{\frac{(k-1)\pi}{n}\leq t \leq \frac{k\pi}{n}}|\phi(t)-\phi_{k}^{*}(t)|\leq5\sqrt{2}\pi\frac{M}{n}\|r\|_{C}, \ n\in\mathbb{N}.
 \end{equation}

 It follows from  (\ref{b1}), (\ref{b21}) and (\ref{b10}) that
 \begin{equation}\label{b11}
 \|\phi\|_{\infty}= \|r\|_{C}\Big(1+\delta_{\infty,n}^{(1)}\frac{M}{n}\Big), \ \ |\delta_{\infty,n}^{(1)}|\leq5\sqrt{2}\pi.
  \end{equation}

Let us prove inequality  (\ref{a2}) for $s=\infty$. By using the inequality  (\ref{a16}) for $s=\infty$, by applying Lagrange theorem, formulas
 (\ref{a70}) and (\ref{b11}),
 we obtain
$$
\frac{1}{2}\|\phi\big(t+\frac{\pi}{n}\big)-\phi(t)\|_{\infty}
\geq
$$
$$
\geq\|\phi\|_{\infty}-\frac{1}{\sqrt{2}}\bigg(\Big\| r\Big(t+\frac{\pi}{n}\Big)-r(t)\Big \|_{\infty}
+3\pi\|r\|_{C}\frac{M}{n}\bigg)\geq
$$
$$
\geq\|\phi\|_{\infty}-\frac{1}{\sqrt{2}}\Big(\frac{\pi}{n}\sup\limits_{t\in\mathbb{R}}\Big|\frac{r'(t)}{r(t)}\Big|\|r\|_{C}  +3\pi\|r\|_{C}\frac{M}{n}\Big)>
\|r\|_{C}\Big(1-\frac{15\pi}{\sqrt{2}}\frac{M}{n}\Big).
$$
Lemma $1$ is proved.

{\bf Remark 1.} {\it In proving of lemma 1 we established  more exact, than   (\ref{delta}) estimates of quantities  $\delta_{s,n}^{(i)}, i=\overline{1,3}$. Namely, we showed that for $ {n\geq{\bigg\{\begin{array}{cc}
4\pi s M, & 1\leq s<\infty, \\
1, \ \ \ \ \ \  \ & s=\infty, \
  \end{array} }
  }$ the following estimates hold }
$$
|\delta_{s,n}^{(1)}|\leq{\left\{\begin{array}{cc}
 \sqrt{2}(5\pi+2)+4, \ \ & 1\leq s<\infty,\\
  5\sqrt{2}\pi,  & \ \ s=\infty,
\end{array} \right.}
$$
$$
-\frac{15\pi+6}{\sqrt{2}}-4 \leq \delta_{s,n}^{(i)}\leq \sqrt{2}(5\pi+2)+4, \ \ \ i=2,3, \ 1\leq s<\infty,
$$
$$
-\frac{15\pi}{\sqrt{2}}\leq \delta_{s,n}^{(i)} \leq 5\sqrt{2}\pi, \ \ \ i=2,3, \  s=\infty.
$$

E-mail: \href{mailto:serdyuk@imath.kiev.ua}{serdyuk@imath.kiev.ua},
\href{mailto:tania_stepaniuk@ukr.net}{tania$_{-}$stepaniuk@ukr.net}


\newpage
\begin{thebibliography}{10}

\bibitem{Stepanets1}
{\sc Stepanets, A.I.:}
Methods of Approximation Theory. VSP: Leiden, Boston  (2005).


\bibitem{Stepanets_Serdyuk_Shydlich}
 {\sc Stepanets', A.I., Serdyuk, A.S., Shidlich, A.L.:} On some new criteria for infinite differentiability of periodic functions. (Ukrainian, English)
 Ukr. Mat. Zh. \textbf{ 59}(10), 1399-1409 (2007); translation in Ukr. Math. J. \textbf{59}(10), 1569-1580 (2007).




\bibitem{Stepanets_Serdyuk_Shydlich2009}
{\sc Stepanets, A.I., Serdyuk, A.S., Shidlich, A.L.:} On relationship between classes of $(\psi, \overline{\beta})$--differentiable functions and Gevrey classes. (Russian, English)
  Ukr. Mat. Zh. \textbf{61}(1), 140-144 (2009); translation in Ukr. Math. J. \textbf{61}(1), 171-177 (2009).

\bibitem{Kuchpel}
{\sc Stepanets, A.I., Kushpel’, A.K.:}
  Convergence rate of Fourier series and best approximations in the space Lp. (English. Russian original)
Ukr. Math. J. \textbf{39}(4), 389-398 (1987); translation from Ukr. Mat. Zh. \textbf{39}(4), 483-492 (1987).





\bibitem{Kuchpel1989}
{\sc  Kushpel’, A.K.:}
Estimates of the widths of classes of analytic functions. (English. Russian original)
Ukr. Math. J. \textbf{41}(4), 493--496 (1989); translation from Ukr. Mat. Zh. \textbf{41}(4), 567--570 (1989).



\bibitem{Teljakovsky 1989}
{\sc Telyakovskii, S.A.:}
Approximation of functions of high smoothness by Fourier sums. (English. Russian original)
Ukr. Math. J. \textbf{41}(4), 444-451 (1989); translation from Ukr. Mat. Zh. \textbf{41}(4), 510-518 (1989).


\bibitem{Temlyakov1990MZ}
  {\sc Temlyakov, V.N.:}
To the question on estimates of the diameters of classes of infinite- differentiable functions. (Russian)
Mat. Zametki \textbf{47}(5), 155-157 (1990).


\bibitem{Temlyakov1990Vekya} 
 {\sc Temlyakov, V.N.:} On estimates of the diameters of classes of infinite- differentiable functions. (Russian) Dokl. razsh. zased. semin. Inst. Prykl. Mat. im.  I.N.Vekua  \textbf{5}(2), 111-114
    (1990).



\bibitem{serdyuk2004zbirnyk}
{\sc Serdyuk, A.S.:}
On one linear method of approximation of periodic functions. (Ukrainian)
Zb. Pr. Inst. Mat. NAN Ukr. \textbf{1}(1), 294--336 (2004).



\bibitem{Serdyuk2005}
{\sc Serdyuk, A.S.:} Approximation of classes of analytic functions by Fourier sums in uniform metric. (Ukrainian, English)
Ukr. Mat. Zh. \textbf{57}(8), 1079-1096 (2005); translation in Ukr. Math. J. \textbf{57}(8), 1275-1296 (2005).


\bibitem{S_S}
{\sc Serdyuk, A.S., Stepanyuk, T.A.:}
Order estimates for the best approximation and approximation by Fourier sums of classes of infinitely differentiable functions. (Ukrainian. English summary)
Zb. Pr. Inst. Mat. NAN Ukr. \textbf{10}(1), 255-282 (2013).

\bibitem{S_S2}
{\sc Serdyuk, A.S., Stepanyuk, T.A.:} Estimates for the best approximations of the classes of infinitely differentiable functions in uniform and integral metrics. (Ukrainian, English)  Ukr. Mat. Zh., \textbf{66}(9), 1244--1256 (2014); English translation in Ukr. Math. J., \textbf{66}(9), 1393--1407 (2015).





\bibitem{Kol}
 Kolmogoroff, A.: Zur Gr\"{o}ssennordnung des Restgliedes
Fourierschen Reihen differenzierbarer Funktionen. (German) Ann. Math.(2),
\textbf{36}(2), 521--526 (1935).



\bibitem{Nikolsky 1946}
{\sc Nikol’skii, S.M.:}
Approximation of functions in the mean by trigonometrical polynomials. (Russian. English summary)
Izv. Akad. Nauk SSSR, Ser. Mat. \textbf{10}, 207-256 (1946).


\bibitem{Teljakovsky1968}
{\sc Telyakovskii, S.A.:} Approximation of differentiable functions by partial sums of their Fourier series. (Russian, English) Mat. Zametki, \textbf{4}(3), 291–300 (1968);
English translation: Mathematical Notes, \textbf{4}(3), 668–673 (1968).

\bibitem{Serdyuk_grabova}
{\sc Hrabova, U.Z., Serdyuk, A.S.:} Order estimates for the best approximations and approximations by Fourier sums of the classes of
$(\psi,\beta)$-differential functions. (Ukrainian, English) Ukr. Mat. Zh., \textbf{65}(9), 1186–1197 (2013); English translation: Ukr. Math. J., \textbf{65}(9), 1319–1331
(2014).


\bibitem{Serdyuk_Stepaniuk2015}
{\sc Serdyuk, A.S., Stepanyuk, T.A.:} Order estimates for the best approximations and approximations by Fourier sums of the classes
of convolutions of periodic functions of low smoothness in the integral metric. Ukr. Mat. Zh., \textbf{66}(12), 1658–1675 (2014); English
translation: Ukr. Math. J., \textbf{65}(12), 1862–1882 (2015).



\bibitem{Stechkin 1980}
{\sc Stechkin, S.B.:}
An estimate of the remainder term of Fourier series for differentiable functions. (Russian)
Tr. Mat. Inst. Steklova 145, 126-151 (1980).










\bibitem{Serdyuk2011}
{\sc Serdyuk, A.S.,  Sokolenko, I.V.:} Uniform approximation of the classes of $(\psi, \overline{\beta})$--differentiable functions by linear methods. (Ukrainian. English summary)
Zb. Pr. Inst. Mat. NAN Ukr. \textbf{8}(1), 181-189 (2011).



\bibitem{Step monog 1987}
{\sc Stepanets, A.I.:}
Classification and approximation of periodic functions. Rev., updated and transl. by P. V. Malyshev and D. V. Malyshev. (English)
Mathematics and its Applications (Dordrecht). 333. Dordrecht: Kluwer Academic Publishers (1995).





\bibitem{Step 1984}
{\sc Stepanets, A.I.:} Deviation of Fourier sums on classes of infinitely differentiable functions. (English)  Ukr. Mat. J.,  \textbf{36}(6), 567--573
(1984).








\bibitem{Korn}
{\sc Korneichuk N.P.:} Exact Constants in Approximation Theory. Encyclopedia of Mathematics and Its Applications, Vol. 38, Cambridge Univ. Press, Cambridge, New York (1990).

\bibitem{Bari}
{\sc Bari, N.K.:}
A treatise on trigonometric series. Vol. I, II. Authorized translation by M. F. Mullins. (English)
Oxford-London-New York-Paris-Frankfurt: Pergamon Press. XXIII, 553 p.; XIX, 508 p. (1964).


\bibitem{Titmarsh}
{\sc Titchmarsh, E.C.:} Introduction to the Theory of Fourier Integrals (2nd.ed.) Oxford University Press (1948).

\bibitem{St-Ruk-Ch}
{\sc Stepanets, A.I., Rukasov, V.I., Chaichenko, S.O.:}
Approximations by the de la Vallee-Poussin sums.  (Russian)
Pratsi Instytutu Matematyky Natsional’noi Akademii Nauk Ukrainy. Matematyka ta ii Zastosuvannya \textbf{68}. Kyiv: Instytut Matematyky NAN Ukrainy (2007).

\end{thebibliography}
\end{document}